\documentclass{article}
\usepackage{amsmath}[1996/11/01]
\usepackage{amssymb,amsthm,amsfonts,amsxtra}
\usepackage{amscd}

\let\catf=\mathbf

\def\[#1\]{\begin{equation}#1\end{equation}}
\makeatletter
\def\beq{%
   \relax\ifmmode
      \@badmath
   \else
      \ifvmode
         \nointerlineskip
         \makebox[.6\linewidth]%
      \fi
      $$
   \fi
}
\def\eeq{%
   \relax\ifmmode
      \ifinner
         \@badmath
      \else
         $$
      \fi
   \else
      \@badmath
   \fi
   \ignorespaces
}

\def\enddisplaymath{\eeq\global\@ignoretrue}
\makeatother

\newtheorem{thm}{Theorem}
\newtheorem{cor}[thm]{Corollary}
\newtheorem{lem}[thm]{Lemma}
\newtheorem{prop}[thm]{Proposition}
\newtheorem{conj}{Conjecture}
\theoremstyle{definition}
\newtheorem{defn}{Definition}
\theoremstyle{remark}
\newtheorem*{rem}{Remark}
\newtheorem*{eg}{Example}
\newtheorem{rems}{Remark}[thm]

\numberwithin{equation}{section}
\numberwithin{thm}{section}

\newcommand{\C}{\mathbb C}
\newcommand{\R}{\mathbb R}
\newcommand{\RP}{\mathbb{RP}}

\renewcommand{\P}{\mathbb P}
\newcommand{\A}{\mathbb A}
\newcommand{\Z}{\mathbb Z}

\newcommand{\F}{\mathbb F}
\newcommand{\Build}{\catf{Build}}
\DeclareMathOperator{\rt}{root}
\DeclareMathOperator{\child}{child}
\DeclareMathOperator{\cone}{cone}

\DeclareMathOperator{\op}{op}
\DeclareMathOperator{\im}{im}

\DeclareMathOperator{\ori}{or}

\DeclareMathOperator{\Tor}{Tor}

\newcommand{\la}{\langle}
\newcommand{\ra}{\rangle}

\begin{document}

\title{The homology of real subspace arrangements} \author{Eric
  M. Rains\footnote{Department of Mathematics, University of California,
    Davis; presently at Department of Mathematics, California Institute of
    Technology}}

\date{December 4, 2009}
\maketitle

\begin{abstract}
Associated to any subspace arrangement is a ``De Concini-Procesi model'', a
certain smooth compactification of its complement, which in the case of the
braid arrangement produces the Deligne-Mumford compactification of the
moduli space of genus $0$ curves with marked points.  In the present work,
we calculate the integral homology of {\em real} De Concini-Procesi models,
extending earlier work of Etingof, Henriques, Kamnitzer and the author on
the ($2$-adic) integral cohomology of the real locus of the moduli space.
To be precise, we show that the integral homology of a real De
Concini-Procesi model is isomorphic modulo its $2$-torsion to a sum of
cohomology groups of subposets of the intersection lattice of the
arrangement.  As part of the proof, we construct a large family of natural
maps between De Concini-Procesi models (generalizing the operad structure
of moduli space), and determine the induced action on poset cohomology.  In
particular, this determines the ring structure of the cohomology of De
Concini-Procesi models (modulo $2$-torsion).

\end{abstract}

\tableofcontents

\section{Introduction}

The main result of \cite{EHKR} was the determination of the cohomology ring
structure of the real locus $\overline{M_{0,n}}(\R)$ of the moduli space of
stable genus 0 curves with $n$ marked points.  It was shown there that the
cohomology of $\overline{M_{0,n}}(\R)$ could be expressed in terms of the
homology of intervals of a certain poset, namely the poset of partitions of
$\{1,2,\dots,n-1\}$ in which each part is of odd size.  This was, of
course, quite reminiscent of the corresponding result (see, e.g.,
\cite{LehrerGI:1987}) for the cohomology of the configuration space of $n$
distinct points in $\P^1(\C)$, which is expressed in terms of the homology
of the poset of {\em all} partitions of $\{1,2,\dots,n-1\}$.  As the latter
result generalizes to arbitrary subspace arrangements
\cite{GoreskyM/MacPhersonR:1988,YuzvinskyS:2002,DeligneP/GoreskyM/MacPhersonR:2000,deLonguevilleM/SchultzCA:2001},
it was natural to look for a corresponding generalization of \cite{EHKR}.

In our context, the appearance of subspace arrangements comes from an
interpretation of $\overline{M_{0,n}}$ (real or complex) as a special case
of a construction of De Concini and Procesi \cite{DeConciniC/ProcesiC:1995}
of a ``wonderful'' compactification associated to an arbitrary subspace
arrangement (in this case, the braid arrangement of type $A_{n-2}$, as
remarked in \cite[p. 483]{DeConciniC/ProcesiC:1995}).  When the spaces in
the arrangement are real, the De Concini-Procesi model is a smooth real
projective variety, and thus gives rise to a smooth manifold.  In the
present note, we determine the (integral) homology groups of these
manifolds.  More precisely, the homology of these manifolds consists of a
large amount of $2$-torsion (analogous to the homology of the spaces
$\RP^n$), which is implicitly determined (via the mod $2$ homology) in
Section \ref{sec:Z2} below, but if we quotient out by this $2$-torsion, the
result can then be expressed (canonically) as a direct sum of cohomology
groups of certain simplicial complexes described in Section
\ref{sec:partitions}.  These simplicial complexes can in turn be subdivided
to obtain the order complexes of certain subposets of the lattice of
subspaces in the arrangement; the primary constraint ($2$-divisibility) on
the subspaces being that they decompose as transverse intersections of
even-codimensional subspaces.  The precise statement of the main theorem is
given as Theorem \ref{thm:main} below.

As might be guessed from the fact that we switched to considering homology
rather than cohomology as in \cite{EHKR}, our techniques are somewhat
different; in particular, we do not have the luxury of an explicit
presentation for the cohomology algebra.  (Similarly, while \cite{EHKR}
made some use of the fact that $\overline{M_{0,n}}(\R)$ is a $K(\pi,1)$,
this fails to hold for arbitrary real De Concini-Procesi models (e.g.,
$\RP^n$).) The approach of \cite{EHKR} of obtaining information about
$2$-adic (co)homology using the mod 2 (co)homology does play a significant
role, however (see Section \ref{sec:Z2} below).  Moreover, one of the main
tools of \cite{EHKR} was the fact that the moduli spaces form an operad,
giving rise to a large collection of natural maps that were used there to
distinguish cohomology classes; one of our main tools (even more so than in
\cite{EHKR}) is the observation that these operad maps can be defined for
general De Concini-Procesi models.  In particular, this includes diagonal
morphisms; thus, although we work exclusively with homology, we still
effectively (if somewhat implicitly) determine the ring structure on
cohomology.

Our primary remaining tool, which was not used in \cite{EHKR}, is a certain
long exact sequence associated to blow-ups of real varieties.  Since every
De Concini-Procesi model is a blow-up of a simpler De Concini-Procesi
model, this allows us to reduce a significant portion of the main theorem
(localized away from the prime 2) to the (trivial) case of products of
projective spaces.  In particular, this allows us to resolve a conjecture
of \cite{EHKR} by showing that the (co)homology of $\overline{M_{0,n}}(\R)$
has no odd torsion.  This holds for arbitrary Coxeter arrangements, see
\cite{arrange_cox}, but again not for general De Concini-Procesi models.
(In fact, even if one restricts one's attention to hyperplane arrangements,
one can arrange for the homology of an arbitrary finite simplicial complex
to appear as a graded piece of the homology of some De Concini-Procesi
model; see the remark following the statement of the main theorem.)

The plan of the paper is as follows.  In Section \ref{sec:models}, we
recall De Concini and Procesi's construction, and introduce the associated
``operad'' maps (together with the main ``composition law'' that they
satisfy); we also show how these maps give rise to a natural grading of the
homology groups by the lattice of subspaces in the arrangement.  In Section
\ref{sec:partitions}, we introduce the corresponding combinatorial data, in
particular allowing us to state the main theorem (Theorem \ref{thm:main}).
The final set-up section, Section \ref{sec:blowup}, establishes the long
exact sequence associated to a real blow-up.  Corollary \ref{cor:DCP_les}
of this section discusses the special case of the blow-up long exact
sequence associated to De Concini-Procesi models, and in particular shows
how this sequence interacts with the natural grading by subspaces.

The proof of the main theorem spans Sections \ref{sec:cells}, \ref{sec:Z2},
and \ref{sec:operadic}.  In Section \ref{sec:cells}, we construct a family
of natural cell structures on the De Concini-Procesi model, which allow us
to construct the isomorphism of the main theorem via a chain map.  It is
then fairly straightforward to show that the chain map respects the blow-up
long exact sequence, and thus by induction that it induces an isomorphism
on homology with coefficients in $\Z[1/2]$.  This argument fails to control
the $2^k$-torsion, however, which is the subject of Section \ref{sec:Z2};
there, we show that the chain map gives an isomorphism (modulo $2$-torsion)
over $\Z/4\Z$, which then implies the isomorphism property over the
$2$-adic integers.  In the process, we give an explicit basis of the mod
$2$ homology consistent with the natural grading (its primary distinction
from the cohomology basis of \cite{YuzvinskyS:1997}), which together with
the main theorem determines the $2$-torsion of the homology groups.  The
proof of the main theorem is finished in Section \ref{sec:operadic}, where
we show that the action of the operad maps on homology agrees with the
action described combinatorially in Section \ref{sec:partitions}.  This in
particular determines the ring structure on cohomology.

Finally, in Section \ref{sec:further}, we discuss some possibilities for
further generalizations; see especially Conjecture \ref{conj:complex},
which considers the case of $\R$-rational (i.e., closed under conjugation)
arrangements.  Also of some interest is Theorem \ref{thm:twist}, which
gives an interpretation of the cohomology of the {\em full} poset of
``$2$-divisible'' subspaces (when that poset has no maximal element) in
terms of the homology, twisted by a certain sheaf, of the De
Concini-Procesi model.

Notational convention: Unless otherwise stated, all modules are
supermodules, with corresponding conventions for tensor products.

Acknowledgements.  The author would like to thank his coauthors P. Etingof,
A. Henriques, and J. Kamnitzer on \cite{EHKR} for introducing him to these
questions, and especially Henriques for discussions relating to blow-ups in
$\Z[1/2]$-cohomology.  In addition, the author would like to thank
S. Devadoss and especially S. Yuzvinsky for motivating discussions, as well
as several referees for useful comments.  The author was supported in part
by NSF Grant No. DMS-0401387.

\section{De Concini-Procesi models}\label{sec:models}

We will need to recall (and extend) some notions from
\cite{DeConciniC/ProcesiC:1995}.  Note that our definition is slightly more
general than theirs, but only in that their construction is not closed
under taking products (see Lemma \ref{lem:product} and surrounding
remarks).  The considerations of functoriality and operadicity are new, and
in particular the construction of a natural grading by subspaces on the
(co)homology groups (see the remark following Corollary \ref{cor:idems}).

A {\em subspace arrangement} in a finite-dimensional vector space $V$ is a
finite collection ${\cal G}$ of subspaces of $V^*$.  Note that this induces
by duality a corresponding collection of subspaces of $V$, but it will be
convenient to use the dual notation.  Given a subspace arrangement ${\cal
  G}$, we define ${\cal C}_{\cal G}$ to be the lattice generated by ${\cal
  G}$; that is, the set of all sums of subsets of ${\cal G}$.  (By
convention, this includes the empty sum; i.e., $0\in {\cal C}_{\cal G}$.)
Note that ${\cal C}_{\cal G}$ is indeed a lattice with respect to
inclusion, but the meet operation in this lattice is not simply the
intersection of subspaces.

Given $U\in {\cal C}_{\cal G}$, an ${\cal G}$-{\em decomposition} of $U$ is
a collection of nonempty subspaces $U_i\in {\cal C}_{\cal G}$ such that
\[
U = \bigoplus_{1\le i\le k} U_i
\]
and for every $G\in {\cal G}$ such that $G\subset U$, $G\subset U_i$ for
exactly one $i$.  Note that every element of ${\cal G}$ is ${\cal
  G}$-indecomposable, and the notion of decomposition depends only on the
collection of ${\cal G}$-indecomposable subspaces.  In particular, if we
let $\overline{\cal G}$ denote that collection, then ${\cal C}_{\cal
  G}={\cal C}_{\overline{\cal G}}$ and an element is ${\cal
  G}$-indecomposable iff it is $\overline{\cal G}$-indecomposable.  A {\em
  building set} is a subspace arrangement ${\cal G}$ such that ${\cal
  G}=\overline{\cal G}$; we have thus seen that every arrangement induces a
building set.

If ${\cal G}$ is a building set in $V$, and ${\cal G}'$ is a building set
in $V'$, a {\em morphism} from ${\cal G}$ to ${\cal G}'$ is a linear
transformation $f:V\to V'$ such that $f^*(G)\in {\cal G}$ for all $G\in
{\cal G}'$.  We thus obtain for every field $K$ a category $\Build(K)$ of
building sets of vector spaces over $K$.  There are two important special
cases of morphisms.  First, if ${\cal G}\supset{\cal G}'$ are both building
sets in $V$, then the identity map on $V$ induces a morphism $\iota:{\cal
  G}\to {\cal G}'$.  Second, we have morphisms $f$ such that
\begin{align}
{\cal G}&=\overline{\{f^*(G):G\in {\cal G}'\}}\label{eq:induced1}\\
{\cal G'}&=f_*({\cal G}):=\{C:C \in {\cal C}_{{\cal G}'}|f^*(C)\in {\cal
  G}\}.\label{eq:induced2}
\end{align}
Since for a general morphism,
\begin{align}
{\cal G}&\supset \overline{\{f^*(G):G\in {\cal G}'\}},\\
{\cal G}'&\subset f_*({\cal G}),\\
\overline{\{f^*(G):G\in {\cal G}'\}} &= \overline{\{f^*(G):G\in f_*({\cal G})\}},
\end{align}
it follows that any morphism can be decomposed as $\iota\circ f\circ
\iota'$ with $f$ satisfying \eqref{eq:induced1},\eqref{eq:induced2}.  Both
types of morphisms can be further decomposed; indeed, it suffices to take
morphisms with $\dim(\ker f)=1$ and ${\cal G}'=f_*({\cal G})$, and
morphisms $\iota$ such that $|{\cal G}'|=|{\cal G}|-1$.  The only
nontrivial thing to prove is the following.

\begin{lem}\label{lem:build_remove_minimal}
Given any two building sets ${\cal G}\supset {\cal G}'$ on the same space,
there exists an element $G\in {\cal G}\setminus {\cal G}'$ such that ${\cal
  G}\setminus \{G\}$ is a building set.
\end{lem}

\begin{proof}
Let $G$ be a minimal element of ${\cal G}$ not in ${\cal G}'$, and follow
the proof of \cite[Prop. 2.5(1)]{DeConciniC/ProcesiC:1995}.
\end{proof}

De Concini and Procesi associate a smooth projective variety to a
building set ${\cal G}$ as follows.  Let ${\cal A}_{\cal G}$ be the affine
variety
\[
{\cal A}_{\cal G} = V \setminus \bigcup_{G\in {\cal G}} G^\perp.
\]
Then for each $G$ we have a map ${\cal A}_{\cal G}\to \P_G$, where
$\P_G=\P(V/G^\perp)$, and we thus have a map
\[
{\cal A}_{\cal G} \to \prod_{G\in {\cal G}} \P_G.
\]
Then ${\overline Y}_{\cal G}$ is the closure of the image of this map.  A
more local description can be given as follows: Let $\rho_G:{\overline
  Y}_{\cal G}\to \P_G$ be the natural map.  Then ${\overline Y}_{\cal
  G}$ is the locus of points $x\in \prod_{G\in {\cal G}} \P_G$ such that
for every pair $H\subset G\in {\cal G}$, $(\rho_G(x),\rho_H(x))$ is in the
(closed) graph of the projection $\P_G\to \P_H$.

If $\dim(G)=1$ (i.e., if $G^\perp$ is a hyperplane), then $\P_G$ is a
single point, so we find that removing all hyperplanes from a building set
has no effect on the corresponding variety.  Conversely, adjoining a
hyperplane to a building set has no effect, so long as the result is still
a building set.  A useful criterion for this is the following.

\begin{lem}\label{lem:hyperplane}
  Suppose $v\in V^*$ is a nonzero vector such that for any $C\in {\cal
    C}_{\cal G}$ containing $v$ there exists $G\in {\cal G}$ such that
  $v\in G\subset C$.  Then ${\cal G}\cup\{\la v\ra\}$ is a building set.
\end{lem}

\begin{proof}
Let ${\cal G}':={\cal G}\cup\{\la v\ra\}$, and consider $C\in {\cal
  C}_{{\cal G}'}$ which is ${\cal G}'$-indecomposable.  We need to show
that $C\in {\cal G}'$.

If $v\notin C$, then $C\in {\cal C}_{\cal G}$.  A ${\cal G}$-decomposition
of $C$ would then also decompose $C$ as an element of ${\cal C}_{{\cal
    G}'}$; it follows that $C\in {\cal G}$.

If $v\in C$, but still $C\in {\cal C}_{\cal G}$, then by assumption, there
exists $G\in {\cal G}$ with $v\in G\subset C$.  In particular, $G$ is
contained in some component of the ${\cal G}$-decomposition of $C$, and
thus $v$ is contained in that component.  It follows that the ${\cal
  G}$-decomposition of $C$ is the same as its ${\cal G}'$-decomposition,
and thus $C\in {\cal G}$.

Finally, we have the case $v\in C$, $C\notin {\cal C}_{\cal G}$.  In
particular, we can write $C=D+\la v\ra$ for some $D\in {\cal C}_{\cal G}$
not containing $v$.  We claim that the ${\cal G}'$-decomposition of $C$ is
obtained by adjoining $\la v\ra$ to the ${\cal G}$-decomposition of $D$.
Suppose otherwise: that there exists $G'\in {\cal G}$ such that $G'\subset
C$ but $G'\not\subset D$.  But then $C=G'+D\in {\cal C}_{\cal G}$, a
contradiction.  It follows that $C=\la v\ra$.
\end{proof}

In particular, if $v\in G\in {\cal G}$ is such that if $v\in C\in {\cal
  C}_{\cal G}$, then $G\subset C$, then ${\cal G}\cap \{\la v\ra\}$ is a
building set.  Over an infinite field, there exists such a vector in any
$G\in {\cal G}$, and we may thus repeatedly adjoin hyperplanes to ${\cal
  G}$ without affecting the variety, until each $G\in {\cal G}$ is the span
of $\dim(G)+1$ hyperplanes in the collection, and those hyperplanes are in
general position.  Letting ${\cal H}$ denote this collection of
hyperplanes, we find $\overline{{\cal H}} = {\cal G}\cup {\cal H}$, and
thus ${\overline Y}_{\cal G}\cong {\overline Y}_{\overline{\cal H}}$.  In
other words, every De Concini-Procesi model is isomorphic to the De
Concini-Procesi model of a hyperplane arrangement.  Thus by considering
general building sets, we do not in fact add any more generality than if we
merely considered the hyperplane case; they do, however, form a useful
tool.

We note that the construction of ${\overline Y}_{\cal G}$ attributed to De
Concini and Procesi above is not {\em quite} the construction they give.
To be precise, they also include the map ${\cal A}_{\cal G}\to \P(V)$, or
equivalently in our notation assume that $V^*\in {\cal G}$, in which case
${\cal A}_{\cal G}$ embeds in ${\overline Y}_{\cal G}$; it will be notationally
convenient to allow the slightly more general case.  Our case reduces
easily to the case $V^*\in {\cal G}$, as follows.  Given $W\subset V^*$,
define the restriction ${\cal G}|_W$ by
\[
{\cal G}|_W = \{G\in {\cal G}|G\subset W\},
\]
a building set in $W^*=V/W^\perp$.  Since $(V/W^\perp)/G^\perp=V/G^\perp$
for $G\subset W$, we immediately obtain the following, using the local
description of ${\overline Y}_{\cal G}$.  The {\em root} of a building set
is the maximal element $\rt({\cal G})$ of the lattice ${\cal C}_{\cal G}$.

\begin{lem}\label{lem:product}
For any $W\subset V^*$, there exists a natural map ${\overline Y}_{\cal
  G}\to {\overline Y}_{{\cal G}|_W}$.  If $G_1\oplus G_2\oplus\cdots\oplus
G_k$ is the decomposition of $\rt({\cal G})$, then
\[
{\overline Y}_{\cal G}\to \prod_{1\le i\le k} {\overline Y}_{{\cal
    G}|_{G_i}}
\]
is an isomorphism.  In particular,
\[
\dim{\overline Y}_{\cal G} = \sum_{1\le i\le k} (\dim(G_i)-1) =
\dim(\rt({\cal G}))-k.
\]
\end{lem}

\begin{rem}
In particular, each of the factors contains the appropriate ambient space.
Also, for topological purposes we note that ${\overline Y}_{\cal G}$ in our
notation is homotopic to the variety $Y_{\cal G}$ (the closure with a
factor $V$ added to the map) discussed in \cite{DeConciniC/ProcesiC:1995}.
\end{rem}

\begin{prop}
The construction ${\overline Y}$ defines a functor from $\Build(K)$
to the category of smooth projective varieties over $K$.
\end{prop}

\begin{proof}
Let $f$ be a morphism from ${\cal G}$ to ${\cal G}'$.  To specify the
associated map ${\overline Y}_f:{\overline Y}_{\cal G}\to {\overline
  Y}_{{\cal G}'}$, it suffices to specify $\rho_G\circ {\overline Y}_f$ for
each $G\in {\cal G}'$.  We simply take
\[
\rho_G\circ {\overline Y}_f=\rho_{f^*(G)}\circ \P(f),
\]
where $\P(f):\P(V/(f^*(G))^\perp)\to \P(V'/G^\perp)$ is the natural
morphism, which is well-defined (and injective) since
$(f^*(G))^\perp=f^{-1}(G^\perp)$.  The local conditions are then
straightforward to verify, as is the fact that ${\overline Y}$ respects
composition of morphisms.
\end{proof}

\begin{rem}
The defining maps $\rho_G$ are associated in this way to the morphisms
$\iota:{\cal G}\to \{G\}$.  Similarly, the diagonal map ${\overline
  Y}_{\cal G}\to {\overline Y}_{\cal G}\times {\overline Y}_{\cal G}$ 
is associated to the diagonal map $\Delta:V\to V\oplus V$; more precisely,
the diagonal map is the composition
\[
\begin{CD}
{\overline Y}_{\cal G}
@>{\overline Y}_\Delta>>
{\overline Y}_{{\cal G}\oplus {\cal G}}
@>\sim>>
{\overline Y}_{\cal G}\times
{\overline Y}_{\cal G}.
\end{CD}
\notag
\]
\end{rem}

In fact, ${\overline Y}$ satisfies a more general version of functoriality,
in that it has a sort of ``operadic'' structure.  Let ${\cal G}$ and ${\cal
  G}'$ be building sets in $V$ and $V'$ respectively.  A {\em weak
  morphism} $f:{\cal G}\to {\cal G}'$ is a linear transformation $f:V\to
V'$ such that $f^*(G)\in {\cal G}\cup\{0\}$ for all $G\in G'$.  

\begin{thm}
Given any weak morphism $f:{\cal G}\to {\cal G}'$, there exists a natural
morphism (the {\em operad map})
\[
\phi_f:
{\overline Y}_{{\cal G}'|_{\ker(f^*)}}
\times
{\overline Y}_{\cal G}
\to
{\overline Y}_{{\cal G}'}.
\]
These maps satisfy the {\em composition law}, which states that given any
two weak morphisms $f:{\cal G}\to {\cal G}'$, $g:{\cal G}'\to {\cal G}''$,
the diagram
\[
\begin{CD}
{\overline Y}_{{\cal G}''|_{\ker{g^*}}}
\times
{\overline Y}_{{\cal G}'|_{\ker(f^*)}}
\times
{\overline Y}_{\cal G}
@>{1\times \phi_f}>> 
{\overline Y}_{{\cal G}''|_{\ker{g^*}}}
\times
{\overline Y}_{{\cal G}'}
\\
@V{\phi_{g|_{\im(f)}}\times 1}VV
@VV{\phi_g}V
\\
{\overline Y}_{{\cal G}''|_{\ker{f^*\circ g^*}}}
\times
{\overline Y}_{\cal G}
@>{\phi_{g\circ f}}>>
{\overline Y}_{{\cal G}''}
\end{CD}
\]
commutes, where $g|_{\im(f)}$ is the induced weak morphism
\[
g|_{\im(f)}:{\cal G}'|_{\ker(f^*)}\to {\cal G}''|_{\ker(f^*\circ g^*)}.
\]
\end{thm}

\begin{proof}
We need to specify $\rho_G\circ \phi_f$ for each $G\in {\cal G}'$.  If
$G\subset \ker(f^*)$, then we simply set $\rho_G\circ\phi_f=\rho_G$,
projecting from ${\overline Y}_{{\cal G}'|_{\ker(f^*)}}$.  Otherwise, we
compose $\rho_{f^*(G)}$ (projecting from ${\overline Y}_{\cal G}$) with the
induced map $\P(f):\P(V/(f^*(G))^\perp)\to \P(V'/G^\perp)$ as before.
\end{proof}

\begin{rems}
  Note that if $f$ is a linear transformation such that $\ker(f^*)\in {\cal
    C}_{{\cal G}'}$, then $\{f^*(G):G\in {\cal G}'|f^*(G)\ne 0\}$ is a
  building set in $V$, called the {\em induced} building set, and denoted
  $f^*({\cal G}')$.  We will denote the corresponding weak morphism from
  $f^*({\cal G}')$ to ${\cal G}'$ by $\tau(f)$, and say that such a weak
  morphism is {\em purely operadic}.  Note that the corresponding operad
  map is injective.  For $C\in {\cal G}'$, we denote by $\phi_C$ the operad
  map associated to $\tau(i_C)$ where $i_C:C^\perp\to V'$ is the inclusion
  map.
\end{rems}

\begin{rems}
In general, any weak morphism can be factored as a product of a morphism
and a purely operadic weak morphism.  Indeed, given a general weak morphism
$f:{\cal G}\to {\cal G}'$, let $C=\rt({\cal G}'|_{\ker f^*})$.  Then we may
factor the linear transformation $f$ as $i_C\circ g$.  But then as a weak
morphism, $f=\tau(i_C)\circ g$.  The composition law in this case simplifies,
and thus we find
\[
\phi_f
=
\phi_C\circ(1\times \overline{Y}_g)
\]
as one might expect.
\end{rems}

\begin{rems}
It will be helpful to note the forms the composition law takes when one of
the maps is a morphism.  If $f$ is a morphism, then
\[
\phi_{g\circ f} = \phi_g\circ (1\times \phi_f)
\]
while if $g$ is a morphism,
\[
\phi_g\circ\phi_f = \phi_{g\circ f}\circ(\phi_{g|_{\im(f)}}\times 1).
\]
Of course, if both are morphisms, the composition law is simply
\[
\phi_{g\circ f}=\phi_g\circ \phi_f.
\]
\end{rems}

\begin{eg}
As an example, consider the braid arrangement
\[
A_{n-1}=\{\langle e_i-e_j\rangle:1\le i<j\le n\}.
\]
The indecomposable subspaces in this arrangement are those of the form
\[
\langle e_i-e_j:i,j\in S\rangle
\]
for some subset $S\subset \{1,2,\dots,n\}$ with $|S|\ge 2$, and general
subspaces are associated to partitions of $\{1,2,\dots,n\}$.  The
associated De Concini-Procesi model is thus constructed as a subvariety of
a product of projective spaces of the form
\[
\prod_{S\subset \{1,2,\dots,n\}} \P^{|S|-2};
\]
to be precise, it is the intersection of the graphs of the maps $\P(S)\to
\P(T)$ for $T\subset S$ given by omitting those coordinates not in $T$.

There is a natural morphism from the moduli space $\overline{M_{0,n+1}}$
(with one marked point singled out as special) to this De Concini-Procesi
model, as follows.  For each subset $S\subset \{1,2,\dots,n\}$, we map the
moduli space to $\P(S)$ by first forgetting all points with labels outside
$S$, collapsing components as necessary to preserve stability, then further
collapse all components not containing the special point.  We thus obtain a
copy of $\P^1$ together with a function from $S$ to $\P^1$ that avoids the
special point.  Equivalently, taking the special point to $\infty$, we
obtain a point in the affine space $\A^{|S|}$, with at least two distinct
coordinates.  The further isomorphisms of $\P^1$ quotient this affine space
by the diagonal subspace, then by scalar multiplication, and thus we obtain
a point in $\P^{|S|-2}$ associated to each stable curve.  The local
conditions are clearly satisfied, so this collection of morphisms to
projective space induce a morphism to ${\overline Y}_{\overline{A_{n-1}}}$.
By considering its action on the open set $M_{0,n+1}$, we find that this
morphism is birational; it is also straightforward to verify that it is
bijective, and thus an isomorphism.
 
The moduli spaces (or, rather, their real or complex loci) form a
topological operad in the usual sense, which is a special case of the above
generalized operad structure.  To be precise, for any composition
$\alpha_1,\dots,\alpha_k$ of $n$ with nonzero parts, there is an associated
weak morphism $f_\alpha:A_{k-1}\to A_{n-1}$ which maps $e_1$ to the sum of
the first $\alpha_1$ basis vectors, $e_2$ to the sum of the next $\alpha_2$
basis vectors, etc.  The kernel of $f_\alpha^*$ is then a subspace in the
arrangement, to wit the subspace associated to the set partition
\[
\{1,2,\dots,\alpha_1\},\{\alpha_1+1,\dots,\alpha_1+\alpha_2\},\dots,
\]
and $f^*(\overline{A_{n-1}})=\overline{A_{k-1}}$.  In particular,
$f_\alpha$ is purely operadic, and the associated operad map is the usual
operad map
\[
\overline{M_{0,\alpha_1+1}}\times
\cdots\times \overline{M_{0,\alpha_k+1}}
\times \overline{M_{0,k+1}}
\to
\overline{M_{0,n+1}}
\]
obtained by gluing the special points of the first $k$ curves to
corresponding non-special points of the last curve.  (Similarly, the
``forget a point'' map is associated to a morphism $A_n\to A_{n-1}$.)
Moreover, the usual operad axiom becomes just a special case of the general
composition law.  Note, however, that the {\em cyclic} operad structure of
the moduli space does not seem to be compatible with its interpretation as
a De Concini-Procesi model, as it does not respect the role of the special
point.
\end{eg}

\begin{prop}
Let $f:{\cal G}\to {\cal G}'$ be a weak morphism with ${\cal G}\subset
f^*({\cal G}')$.  Then $\phi_f$ is injective.
\end{prop}

\begin{proof}
Simply observe that if $f^*(G')=G$, then the map $\rho_{G'}$ on the
codomain is the composition of the map $\rho_G$ on the domain with an
injection.
\end{proof}

\begin{prop}
Let $f:{\cal G}\to {\cal G}'$ be a surjective (weak) morphism.  Then
$\phi_f$ is surjective, and is birational iff $\rt({\cal G})=f^*(\rt({\cal
  G}'))$ and both roots have the same number of components.
\end{prop}

\begin{proof}
The image of the restriction
\[
f:{\cal A}_{\cal G}\to {\cal A}_{{\cal G}'}
\]
is dense if $f:V\to V'$ is surjective, and thus the corresponding map of
projective varieties must be surjective; birationality then follows by
comparing dimensions.
\end{proof}

\begin{rem}
In particular, this applies to any morphism of the form $\iota:{\cal G}\to
{\cal G}'$.
\end{rem}

\begin{cor}\label{cor:idems}
Given a building set ${\cal G}$ and a space $C\in {\cal C}_{\cal G}$, the
natural surjection
\[
\pi_C:=\phi_\iota:\overline{Y}_{\cal G}\to \overline{Y}_{{\cal G}|_C}
\]
has a natural homotopy class $\pi_C^{[-1]}$ of splittings.  Moreover,
these maps satisfy the identities (up to homotopy)
\begin{align}
\pi_{C\wedge D}\circ \pi_C&=\pi_{C\wedge D}\\
\pi^{\vphantom{[-1]}}_C\circ \pi_D^{[-1]}
&\sim
\pi^{[-1]}_{C\wedge D}\circ \pi^{\vphantom{[-1]}}_{C\wedge D}\\
\pi^{[-1]}_C\circ \pi^{[-1]}_{C\wedge D}
&\sim
\pi^{[-1]}_{C\wedge D}.
\end{align}
\end{cor}

\begin{proof}
The splitting maps arise from
\[
\phi_C:\overline{Y}_{{\cal G}|_C}\times \overline{Y}_{{\cal G}/C}\to
\overline{Y}_{{\cal G}}
\]
by choosing a point in $\overline{Y}_{{\cal G}/C}$.  

Now, the composition law implies that
\[
\pi_C\circ \phi_D
=
\phi_{C\wedge D}\circ (\pi_{C\wedge D}\times (\phi_f\circ \pi_{(C+D)/D}))
:
\overline{Y}_{{\cal G}|_D}\times \overline{Y}_{{\cal G}/D}
\to
\overline{Y}_{{\cal G}|_C},
\]
where $f$ is the natural morphism
\[
f:({\cal G}/D)|_{(C+D)/D}\to ({\cal G}|_C)/(C\wedge D).
\]
If $C=D$ then the right-hand side is just the projection
\[
{\overline Y}_{{\cal G}|_C}
\times
{\overline Y}_{{\cal G}/C}
\to
{\overline Y}_{{\cal G}|_C}
\]
and we thus obtain the desired splitting.  More generally, if we choose a
point in $\overline{Y}_{{\cal G}/D}$, we obtain the second identity.  The
other two identities follow similarly from the composition law.
\end{proof}

\begin{rem}
In particular, the associated retractions satisfy
\[
(\pi_C^{[-1]}\circ \pi^{\vphantom{[-1]}}_C)\circ (\pi_D^{[-1]}\circ \pi^{\vphantom{[-1]}}_D)
\sim
\pi_{C\wedge D}^{[-1]}\circ\pi^{\vphantom{[-1]}}_{C\wedge D},
\]
and thus commute up to homotopy.  It follows that any associated
(co)homology group is graded by ${\cal C}_{\cal G}$; that is, splits as a
direct sum indexed by ${\cal C}_{\cal G}$.  To be precise, for any $C\in
{\cal C}_{\cal G}$, let
\[
H_*({\overline Y}_{\cal G})[C]
\]
denote the subspace of consisting of homology classes fixed by
$\pi_C^{[-1]}\circ \pi_C$ and annihilated by $\pi_D^{[-1]}\circ \pi_D$ for
$D\subsetneq C$.  (Equivalently, these are the classes in the image of
$\pi_C^{[-1]}$ which are annihilaged by $\pi_D$ for $D\subsetneq C$.)
Since these are retractions, so idempotents on
(co)homology, and commute, we find that
\[
H_*({\overline Y}_{\cal G})
=
\bigoplus_{C\in {\cal C}_{\cal G}} H_*({\overline Y}_{\cal G})[C],
\]
and more generally
\[
(\pi_C^{[-1]}\circ \pi_C)H_*({\overline Y}_{\cal G})
=
\bigoplus_{D\subset C} H_*({\overline Y}_{\cal G})[D].
\]
\end{rem}

\begin{prop}
If $V^*\notin {\cal G}$, then ${\overline Y}_{{\cal G}\cup V^*}$ is a
projective space bundle over ${\overline Y}_{\cal G}$.
\end{prop}

\begin{proof}
Let $G_1,\dots,G_m$ be the maximal elements of ${\cal G}$; note that by
adjoining suitable hyperplanes, we can force $\sum_i G_i = V^*$ without
changing ${\overline Y}_{\cal G}$.  We thus need to show that the set of
points in $\P(V)$ compatible with a given point in ${\overline Y}_{\cal G}$
form a projective space.  But the compatibility condition is simply that
the projection to each $\P_{G_i}$, if defined, has the correct value.  If
we choose representatives $p_i\in V/G_i^\perp \setminus 0$, a point in
$V\setminus 0$ is compatible iff it is of the form $\sum_i \alpha_i p_i$;
it follows that the preimage is $\P^{m-1}$ as desired.
\end{proof}

\begin{prop}\label{prop:blowup}
If ${\cal G}={\cal G}'\cup\{G\}$, and $G$ is not a maximal element of
${\cal G}$, then $\overline{Y}_{\cal G}$ is the blow-up of
$\overline{Y}_{{\cal G}'}$ along the image $d'_G$ of the injective map
\[
\phi_G:\overline{Y}_{{\cal G}'|_G}\times \overline{Y}_{{\cal G}/G}
\to \overline{Y}_{{\cal G}'}.
\]
The exceptional divisor is the image $d_G$ of the injective map
\[
\phi_G: \overline{Y}_{{\cal G}|_G}\times \overline{Y}_{{\cal G}/G}
\to \overline{Y}_{{\cal G}}.
\]
\end{prop}

\begin{proof}
The composition law gives
\[
\phi_\iota \circ \phi_G = \phi_G\circ(\phi_\iota\times 1),
\]
and thus $\phi_\iota(d_G)=d'_G$.  Since $\dim(d_G)=\dim(\overline{Y}_{\cal
  G})-1>\dim(d'_G)$ and $d_G$ is a projective space bundle over $d'_G$, it
remains only to check that $\phi_\iota$ is injective on the complement of
$d_G$, which is immediate.
\end{proof}

\begin{rem}
  As in \cite{DeConciniC/ProcesiC:1995} (which used two special cases of
  this proposition), this immediately gives an inductive proof that
  $\overline{Y}_{\cal G}$ is a smooth, irreducible variety, since by the
  above proposition and Lemma \ref{lem:build_remove_minimal}, it can be
  obtained from a product of projective spaces by a sequence of blow-ups;
  it also follows that $\overline{Y}_{\cal G}(\R)$ is a smooth, connected
  manifold for ${\cal G}\in \Build(\R)$.  Also note the consequence that
  the normal bundle of $d'_G$ is trivial if $G$ is minimal in ${\cal G}$,
  since then $d_G$ is a product bundle.

  Note that this consturction of $\overline{Y}_{\cal G}$ as an iterated
  blowup can almost certainly be generalized.  Indeed, Keel's construction
  \cite{KeelS:1992} of $\overline{M}_{0,n+1}$ as a blowup of $(\P^1)^{n-1}$
  is not of the above form.
\end{rem}

In general, for any category $\catf{C}$, we define a {\em universal operad in
$\catf{C}$} to be a functor $F:\Build\to\catf{C}$ that also associates a
morphism
\[
\phi_f:F({\cal G}'|_{\ker(f^*)}\oplus {\cal G})\to F({\cal G}')
\]
to every weak morphism, satisfying the composition law
\[
\phi_g\circ \phi_{1\oplus f}
=
\phi_{g\circ f}\circ \phi_{g|_{\im(f)}\oplus 1}.
\]
A natural transformation between universal operads will be said to be
{\em operadic} if it is compatible with the composition law in the obvious
way.  A {\em universal cooperad in $\catf{C}$} is simply a universal operad in
$\catf{C}^{\op}$.

\begin{rem}
  A universal operad in a tensor category, equipped with natural
  isomorphisms $F({\cal G}\oplus {\cal G}')\to F({\cal G})\otimes F({\cal
    G}')$ (compatible with symmetry and associativity), induces (by
  restriction to the braid arrangements) an operad in the usual sense; this
  is not true for a general universal operad, however.  Also, given a {\em
    topological} universal operad, we may take the homology of $F$ (with
  appropriate coefficients), and thus obtain a universal operad in the
  appropriate category of modules (or, taking cohomology, a universal
  cooperad of rings, assuming the topological universal operad respects
  products).
\end{rem}

\section{Poset homology and operad maps}\label{sec:partitions}

In this section, and until further notice, we restrict our attention to the
case ${\cal G}\in \Build(\R)$.  In this case, there is a significant
difference between odd-dimensional and even-dimensional elements of ${\cal
  G}$: for $\dim(G)>1$, $\P_G(\R)$ is orientable if and only if $\dim(G)$
is even.  This suggests that even-dimensional elements will have particular
significance in the homology of ${\overline Y}_{\cal G}(\R)$.

More generally, let $\Pi^{(m)}_{\cal G}$ be the subposet of ${\cal C}_{\cal
  G}$ consisting of elements $A$ that can be written as direct sums of
elements $G\in {\cal G}$ with $\dim(G)$ a multiple of $m$.  (In the case of
the braid arrangement ${\cal G}=\overline{A_{n-1}}$, $\Pi^{(m)}_{\cal G}$
is the poset of partitions of $\{1,2,\dots,n\}$ into parts all of size
congruent to 1 modulo $m$.)  For the remainder of this section, we fix a
choice of $m$; after this section, we will in fact take $m=2$ almost
exclusively.

For any element $A\in \Pi^{(m)}_{\cal G}$ and any commutative ring $R$, we
define $H_*([0,A];R)$ to be the homology of the chain complex
$C_*([0,A];R)$ in which $C_{k+1}([0,A];R)$ is the free $R$-module spanned
by chains
\[
(0<A_1<\cdots<A_k<A)
\]
in $\Pi^{(m)}_{\cal G}$, and the boundary map is defined by
\[
\partial(0<A_1<\cdots<A_k<A)
=
\sum_i (-1)^i
(0<A_1<\cdots<\hat{A_i}<\cdots<A_k<A).
\]
(Aside from a shift in degree, this is the reduced homology of the order
complex of the open interval $(0,A)$.)  By convention, $C_*([0,0];R)=0$
except in degree 0, where it is spanned by the single chain $(0)$.

It will be convenient to consider an alternate chain complex giving the
same homology, defined in terms of ``${\cal G}$-forests'' (these are called
``nested sets'' in \cite{DeConciniC/ProcesiC:1995}, but we feel ``forest''
is more evocative of the relevant combinatorics, especially in the braid
case).

\begin{defn}
A ${\cal G}$-forest is a subset $F\subset {\cal G}$ such that every
collection of pairwise incomparable elements of $F$ forms a decomposition;
the {\em root} of $F$ is the space
\[
\rt(F):=\sum_{G\in F} G\in {\cal C}_{\cal G}.
\]
(Note that this subspace will not be an element of $F$ in general, but
merely a direct sum of elements of $F$.)  Given $G\in F$, the {\em child}
of $G$ in $F$ is the space
\[
\child_F(G):=\sum_{\substack{H\in F \\ H\subsetneq G}} H.
\]
$F$ is said to be $m$-{\em divisible} if every element of $F$ has dimension
a multiple of $m$.
\end{defn}

\begin{rem}
There is a natural bijection between $\overline{A_{n-1}}$-forests and
forests in the usual sense in which the leaves are labelled $1$,\dots,$n$.
Given such a forest, if one labels each internal node by the set of its
descendant leaves, the resulting collection of subsets forms a
$\overline{A_{n-1}}$-forest.
\end{rem}

The name ``forest'' is justified by the following lemma.

\begin{lem}
Let $F$ be a ${\cal G}$-forest.  For any element $G\in F$, the set of
elements of $F$ containing $G$ forms a chain.
\end{lem}

\begin{proof}
Indeed, if $H_1$, $H_2\in F$ are incomparable elements, then by definition
they form a decomposition; in particular, $G$ cannot be contained in both.
\end{proof}

In other words, $F$ has a natural structure of a forest with nodes labelled
by elements of ${\cal G}$, compatible with inclusion; the root in our sense
is simply the direct sum of the labels of the roots.

Now, for $A\in \Pi^{(m)}_{\cal G}$, define a chain complex $C^f_*(A)$ as
follows: The $R$-module $C^f_k(A)$ is spanned by {\em ordered}
$m$-divisible forests $F$ with root $A$ and $k$ nodes, but with different
orderings identified, up to the obvious sign factor.  For $G\in {\cal G}$,
we define
\[
\partial_G F = \begin{cases}
(-1)^{i-1} F\setminus G,& G=F_i\\
0,& G\notin F;
\end{cases}
\]
the boundary map is then given by the sum of $\partial_G$ with $G$ ranging
over {\em proper} subspaces of components of $A$.  We will also need a
concatenation operation
\[
G\cdot (F_1,F_2,\dots,F_n) = (G,F_1,F_2,\dots,F_n),
\]
or 0 if the result is not a forest.

The following proof is adapted from \cite[Sec. 2.6]{RobinsonA:2004}, which
essentially considered the case $\Pi^{(1)}_{A_n}$ with $A=V$.

\begin{thm}
For all $A\in \Pi^{(m)}_{\cal G}$, there is a canonical isomorphism
\[
H_*([0,A])\cong H^f_*(A).
\]
\end{thm}

\begin{proof}
If $A=0$, the result is immediate (in both cases, $H_0=\Z$ and all other
homology groups are trivial); we may thus assume $A\ne 0$, and thus both
complexes are trivial in degree 0.

Now, given a chain
\[
0<A_1<A_2<\dots<A_k<A,
\]
consider the (partially closed) simplex of numbers $0\le
\tau_k<\cdots<\tau_1<\tau_0=1$.  The remainder of the closure of this
simplex is naturally identified with the disjoint union of simplices
corresponding to chains with steps removed (if $\tau_k=\tau_{k-1}$, remove
$A_k$); as a result, we can glue together all of the simplices to obtain a
geometric simplicial complex $\Sigma$.  The result is simply the the order
complex of $(0,A]$, so its local homology at the point $(0<A)$ is
\[
H_*(\Sigma,\Sigma\setminus\{(0<A)\}) = H_{*+1}([0,A]).
\]

Similarly, given a forest $F$, consider labellings $\tau'$ of the nodes,
subject to the condition that the labels sum to 1, are nonnegative, and the
labels of all non-roots are positive.  Again, if we include labellings of
subforests (with the same root), with the convention that $\tau'=0$ for
removed nodes, we obtain a closed simplex, and a geometric simplicial
complex $\Sigma^f$.  If $x$ is the centroid of the simplex corresponding to
the forest $A$, the local homology at $x$ is
\[
H_*(\Sigma^f,\Sigma^f\setminus\{x\}) = H^f_{*+1}(A).
\]
(Taking local homology with respect to $x$ is equivalent to restricting the
chain complex to those simplices containing the open simplex that contains
$x$.)

The theorem will follow if we can establish a pointed homeomorphism
$\Sigma\cong\Sigma^f$.  Define a (discontinuous) function $\rho:(0,1)\to
\Pi^{(m)}_{\cal G}$ by $\rho(t)=A_l$ whenever $\tau_l<t<\tau_{l-1}$ (with
the convention $A_{k+1}=A$, $\tau_{k+1}=0$.)  Now, if $F$ is the forest
consisting of all components of the $A_k$, we initially label a node $G$ by
the difference between the lim sup and the lim inf of those $t$ for which
$G$ is a component of $\rho(t)$.  The only way such a label can be 0 is if
$\tau_0=0$ and $G$ is a component of $A$ but not $A_1$; in particular
non-root nodes have a nonzero label.  Since the sum of the labels is
uniformly bounded away from 0 and $\infty$, we may rescale to make the sum
1, and obtain the desired labelling $\tau'$.  To invert, we rescale so that
the highest-weight path from the root has weight 1, and can then
immediately recover the labelled chain.

If we begin with the chain $(0<A)$, we find $\rho\equiv A$.  The
corresponding forest is just the set of components of $A$, each labelled (after
normalization) by $1/l$ (if $A$ has $l$ components).  But this is indeed
the centroid of the simplex corresponding to that forest.
\end{proof}

\begin{rem}
The general case $m=1$ was established by different means in
\cite{FeichtnerEM/MuellerI:2005}.
Similarly, the case $m>1$, ${\cal G}={\overline A_{n-1}}$ was established
in \cite{DelucchiE:2005}.
\end{rem}

Under this isomorphism, each simplex of $\Sigma^f$ is identified with a
union of simplices of $\Sigma$, and thus the above isomorphism induces a
chain map.

\begin{cor}
Define a map $\sigma:C^f_*(A)\to C_*([0,A])$ by
\[
\sigma(F):=\sum_\pi \sigma(\pi)
(0<F_{\pi(1)}<F_{\pi(1)}+F_{\pi(2)}<\cdots<F_{\pi(1)}+\cdots+F_{\pi(k-1)}<A),
\]
where the sum is over permutations of the nodes such that if
$F_{\pi(i)}\subset F_{\pi(j)}$, then $i\le j$.  Then $\sigma$ is a chain map
inducing an isomorphism on homology.
\end{cor}

A chain of the form
\[
(0<F_{\pi(1)}<F_{\pi(1)}+F_{\pi(2)}<\cdots<F_{\pi(1)}+\cdots+F_{\pi(k-1)}<A)
\]
will be called a {\em forest} chain.

\begin{lem}
A linear combination of forest chains in $C_*([0,A])$ is in the image of
$\sigma$ if and only if its boundary is also a linear combination of forest
chains.
\end{lem}

\begin{proof}
This follows from geometric considerations, but can also be shown directly
as follows.  Suppose we order $F$ in such a way that $F_i\subset F_j$
implies $i\le j$, and suppose $F_l$, $F_{l+1}$ are incomparable.  Then the
non-forest chain
\[
(0<F_1<\cdots<F_1+\cdots+F_{l-1}<F_1+\cdots+F_{l+1}<\cdots<A)
\]
occurs in the differential of precisely two forest chains, namely
\[
(0<F_1<\cdots<F_1+\cdots+F_{l-1}<F_1+\cdots+F_{l-1}+F_l<F_1+\cdots+F_{l+1}<\cdots<A)
\]
and
\[
(0<F_1<\cdots<F_1+\cdots+F_{l-1}<F_1+\cdots+F_{l-1}+F_{l+1}<F_1+\cdots+F_{l+1}<\cdots<A).
\]
Thus the coefficient of this non-forest chain in the differential is 0 iff
the coefficients of the forest chains are negatives of each other.

Since any two such orderings of $F$ can be connected by a sequence of such
transpositions, the claim follows.
\end{proof}

\begin{rem}
In particular, a chain map between two complexes $C_*([0,A])$ that takes
forest chains to forest chains pulls back to a chain map on the associated
forest complexes.
\end{rem}

Now, define the Whitney homology
\[
W^{(m)}_*({\cal G})=\bigoplus_{A\in \Pi^{(m)}_{\cal G}} H_*([0,A]),
\]
a module graded by both the degree in homology and by the poset
$\Pi^{(m)}_{\cal G}$.  This, of course, can be computed as the homology
of the corresponding sum of poset or forest complexes; thus define, for
instance, $C^W_*({\cal G})=\bigoplus_{A\in {\cal C}_{\cal G}} C^W_*([0,A])$,
and similarly for forests.

\begin{thm}
The Whitney homology extends to a cooperad such that the maps $\phi_f^*$
are homogeneous and respect the poset grading:
\[
\phi_f^*(W^{(m)}({\cal G}')[A])
=
W^{(m)}({\cal G}'|_{\ker(f^*)}\oplus {\cal G})[(A\cap \ker f^*)\oplus
f^*(A)];
\]
in particular $\phi_f^*$ vanishes unless
\[
(A\cap\ker f^*)\oplus f^*(A)\subset \Pi^{(m)}_{{\cal G}'|_{\ker(f^*)}\oplus
{\cal G}}.
\]
\end{thm}

\begin{proof}
Let $f:{\cal G}\to {\cal G}'$ be a weak morphism.
We define a chain map
\[
\phi_f^*
:
C^W_*({\cal G}')
\to
C^W_*({\cal G}'|_{\ker(f^*)})
\otimes
C^W_*({\cal G})
\]
on chains $0<\cdots<A_i<\cdots<A_k=A$ in $\Pi^{(m)}_{\cal G}$ as follows.
If $A\cap\ker f^*$ is not in the chain, we set $\phi_f^*=0$.  Otherwise,
let $l$ be the index of $A\cap\ker f^*$ in the chain, and define
\begin{align}
\phi_f^*(0<\dots<A_i<\dots<A)
&=
\phantom{{}\otimes{}}(0<\dots<A_i\cap\ker f^*<\dots<A_l\cap\ker f^*=A\cap\ker f^*)\notag\\
&\phantom{{}={}}\otimes
(0=f^*(A_l)<\dots<f^*(A_{l+i})<\dots<f^*(A))\notag\\
&\in
C^W_l({\cal G}'|_{\ker(f^*)})
\otimes
C^W_{k-l}({\cal G}).\notag
\end{align}
Since deleting a step cannot introduce $A\cap\ker f^*$ to the chain, the
differential cannot leave the set of ``bad'' chains (i.e., in the kernel of
$\phi_f^*$).  And the only way to produce a bad chain by deleting a node is
to delete $A\cap\ker f^*$; but the corresponding term is missing from the
differential on the image space.  Therefore $\phi_f^*$ is indeed a chain
map.  Moreover, it is easily seen to satisfy the composition law.  Finally,
we can extend it to a chain map
\[
\phi_f^*
:
C^W_*({\cal G}')
\to
C^W_*({\cal G}'|_{\ker(f^*)}\oplus {\cal G})
\]
by composing with the shuffle product
\[
C^W_*({\cal G}'|_{\ker(f^*)})
\otimes
C^W_*({\cal G})
\to
C^W_*({\cal G}'|_{\ker(f^*)}\oplus {\cal G}),
\]
obtaining the map on homology as required.
\end{proof}

\begin{rems}
Since the Whitney homology is a direct sum, it is equivalent to consider
only the induced maps
\[
\phi^*_f:
H_*([0,A]_{{\cal G}'})
\to
H_*([0,A\cap \ker f^*\oplus f^*(A)]_{{\cal G}'|_{\ker(f^*)}\oplus {\cal G}})
\]
where $A\cap \ker f^*\oplus f^*(A)\in \Pi^{(m)}_{{\cal G}'|_{\ker
    f^*}\oplus {\cal G}}$, and similarly for the chain maps.
\end{rems}

\begin{rems}
  The case of the diagonal morphism $\Delta:{\cal G}\to {\cal G}\oplus
  {\cal G}$ is of particular interest:
\[
\phi^*_{\Delta}(0<\cdots<(A_i,B_i)<\cdots< (A,B))
=
\begin{cases}
0<\cdots<A_i+B_i<\cdots A+B & \text{if $A\cap B=0$}\\
0& \text{otherwise}.
\end{cases}
\]
Composing with the shuffle product induces a ring structure on the Whitney
homology, graded by $\Pi^{(m)}_{\cal G}$.
\end{rems}

Since $\phi^*_f$ takes forest chains to forest chains, we conclude the
following.

\begin{cor}
For each weak morphism $f:{\cal G}\to {\cal G}'$, there is a chain map
\[
\phi^*_f
:
C^{Wf}_*({\cal G}')
\to
C^{Wf}_*({\cal G}'|_{\ker(f^*)}\oplus {\cal G})
\]
producing a commutative diagram
\[
\begin{CD}
C^{Wf}_*({\cal G}')@>\phi^*_f>>C^{Wf}_*({\cal G}'|_{\ker(f^*)}\oplus {\cal G})\\
@V\sigma VV @V\sigma VV\\
C^W_*({\cal G}')@>\phi^*_f>>C^W_*({\cal G}'|_{\ker(f^*)}\oplus {\cal G})
\end{CD}.
\]
\end{cor}

Thus, in addition to the cooperad structure on the Whitney homology itself,
we also obtain two cooperads of chain complexes, and an operadic homotopy
between them.  It is unclear, however, how to define the cooperad structure
on the forest complex without passing through the poset complex.

It turns out that the relation to De Concini-Procesi models is more
convenient in terms of poset and forest {\em co}homology.  Of course, the
cooperads of homology chain complexes immediately dualize to operads on the
cohomology chain complexes, and thus give rise to an operad on the Whitney
cohomology (defined in the obvious way).

Even this is not quite the right structure, however.  Note that the induced
chain maps
\[
(\phi_f)_*
:
C^*([0,B\oplus C]_{{\cal G}'|_{\ker(f^*)}\oplus {\cal G}})
\to
C^*([0,A]_{{\cal G}'})
\]
are 0 unless $B=A\cap \ker f^*$, $C=f^*(A)$; in other words, we must have
a short exact sequence
\[
\begin{CD}
0@>>>B@>>> A@>f^*>> C@>>>0.
\end{CD}
\]
In particular, $\dim(A)=\dim(B)+\dim(C)$, and we may thus shift the degrees
of the complex by this dimension without affecting homogeneity of
$(\phi_f)^*$.  To be precise, we will use the chain complex
\[
C^{\dim A-*}([0,A]_{{\cal G}'}),
\]
which becomes a homology complex since the differential now decreases the
degree.

A more subtle correction is a certain twisting of the operad structure.
For any real vector space $V$, let $\ori(V):=H_{\dim(V^*)}(V^*,V^*\setminus
\{0\})$ be the corresponding orientation module; note in particular the
canonical isomorphisms
\[
\ori(V)\cong \tilde{H}_{\dim(V)-1}(S(V^*)),
\]
and when $\dim(V)$ is even,
\[
\ori(V)\cong \tilde{H}_{\dim(V)-1}(\P(V^*)).
\]
In any event, every short exact sequence $0\to V\to W\to X\to 0$ induces a
canonical isomorphism $\ori(V)\otimes \ori(X)\to\ori(W)$.  Thus every
nonzero $\phi^*_f$ induces an isomorphism $\ori(A)\cong \ori(B\oplus C)$,
and we obtain an operad map
\[
(\phi_f)_*
:
C^{\dim(B\oplus C)-*}
([0,B\oplus C]_{{\cal G}'|_{\ker(f^*)}\oplus {\cal G}})\otimes
\ori(B\oplus C)
\to
C^{\dim(A)-*}([0,A]_{{\cal G}'})\otimes \ori(A),
\]
as well as associated maps on (degree-reversed) cohomology.

We may now state our main theorem.  Recall that the notation $[A]$ refers
to the $A$-graded piece of the homology; see the remark following Corollary
\ref{cor:idems}.  Also, we identify the real algebraic variety ${\overline
  Y}_{\cal G}$ with its real locus ${\overline Y}_{\cal G}(\R)$, whenever
this can be done without causing confusion.

\begin{thm}\label{thm:main}
  Let ${\cal G}$ be a real building set, and $A\in \Pi_{\cal G}$.  If
  $A\notin \Pi^{(2)}_{\cal G}$, then
\[
2 H_*({\overline Y}_{\cal G})[A]=0;
\]
otherwise, there is a natural, operadic, isomorphism
\[
H^{\dim(A)-*}_{(2)}([0,A]_{\cal G})\otimes \ori(A)
\to
2 H_*({\overline Y}_{\cal G})[A].
\]
\end{thm}

\begin{rems}
  A similar result (minus the operad structure, and requiring some
  additional hypotheses for naturality) was already known
  \cite{GoreskyM/MacPhersonR:1988} (see also
  \cite{YuzvinskyS:2002,DeligneP/GoreskyM/MacPhersonR:2000,deLonguevilleM/SchultzCA:2001}
  for the ring structure of cohomology) for the complement of a real
  subspace arrangement, namely an isomorphism of its homology with
\[
\bigoplus_{A\in \Pi^{(1)}_{\cal G}} H^{\dim(A)-*}_{(1)}([0,A]_{\cal
  G})\otimes \ori(A).
\]
One curious consequence of this similarity is that if ${\cal G}$ is
obtained from a complex building set by viewing each space as a real space
of twice the dimension, then (since every subspace now has even dimension)
there is an isomorphism
\[
2 H_*({\overline Y}_{\cal G}) \cong H_*(V^*-{\cal G}).
\]
In particular, the homology of ${\overline Y}_{\cal G}(\R)$ can be
arbitrarily complicated (even for hyperplane arrangements, by the remark
following Lemma \ref{lem:hyperplane}), since the same holds for the
complements of complex subspace arrangements.  Indeed, any finite
simplicial complex is homeomorphic to the order complex of a finite atomic
lattice (its lattice of faces), and thus, reversing inequalities, to the
order complex of a finite coatomic lattice.  But any finite coatomic
lattice can be represented as a lattice of subspaces of the space of
complex-valued functions on the coatoms (each element corresponds to the
space of functions vanishing on the coatoms bounding it).
\end{rems}

\begin{rems}
  Dually, there is a natural isomorphism from the cohomology (modulo
  $2$-torsion) of ${\overline Y}_{\cal G}(\R)$ to the (suitably twisted)
  Whitney homology of $\Pi^{(2)}_{\cal G}$, and operadicity (applied to the
  diagonal map) implies that this is an isomorphism of (poset-graded) rings.
\end{rems}

The main theorem has an important consequence for the moduli space of
stable genus 0 curves.

\begin{cor} \cite[Conj. 2.13]{EHKR} The groups $2
  H_*(\overline{M_{0,n}}(\R);\Z)$ and $2 H^*(\overline{M_{0,n}}(\R);\Z)$
  are free.
\end{cor}

\begin{proof}
  As remarked above, we have an isomorphism $\overline{M_{0,n+1}}\cong
  {\overline Y}_{\overline{A_{n-1}}}$ of algebraic varieties, and thus of
  their respective real loci.  It follows that $2
  H_*(\overline{M_{0,n}}(\R),\Z)$ is isomorphic to a direct sum of
  cohomology groups of subposets of the poset of set partitions with all
  parts odd.  But this poset is Cohen-Macaulay, and thus its cohomology is
  free (and supported in the appropriate degree).  The remaining claim
  follows from the universal coefficient theorem.
\end{proof}

\begin{rems}
  The same argument applies to an arbitrary Coxeter arrangement, see
  \cite{arrange_cox}, again because the relevant posets are Cohen-Macaulay.
\end{rems}

\begin{rems}
Note, however, that the explicit presentation of the ring
\[
H^*(\overline{M_{0,n}}(\R);\Z)/H^*(\overline{M_{0,n}}(\R);\Z)[2]
\]
given in \cite{EHKR} does not follow from the present methods.  It is, of
course, far too much to hope for such a presentation for completely general
building sets (or even for general hyperplane arrangements, by the
discussion following Lemma \ref{lem:hyperplane}); for instance, the
cohomology ring need not be generated in degree 1 in that case.
\end{rems}

\section{Blow-ups and homology}\label{sec:blowup}

The construction of ${\overline Y}_{\cal G}$ via repeated blow-ups turns
out to have extremely useful consequences in homology.  As we are
interested in the topological consequences of this, it will be useful to
have a more topological (or, more precisely, differentiable rather than
algebraic) version of blowing-up.  In addition to real and complex
blow-ups, corresponding to algebraic blowing-up of real or complex
varieties, there is a third, ``spherical'' blow-up that it will be useful
to consider.  (For an early application of such blow-ups to subspace
arrangements, see \cite[\S 5]{AxelrodS/SingerIM:1994}.)  In fact, the
spherical blow-up is in some sense universal; the other blow-ups can be
constructed as quotients of the spherical blow-up.

The basic idea of the spherical blow-up is to replace a submanifold $Y$ by
the sphere bundle $N_Y(X)/\R^+$.  That is, the blowup is a new space
$\tilde{X}$ equipped with a continuous (smooth) map to $X$ which is an
isomorphism outside $Y$, and such that the preimage of a given point in $Y$
is identified with the space of unit normal vectors to $Y$ at that point.
The result is homotopic to the complement of $Y$ in $X$, but has the merit
of being compact and {\em almost} smooth; the failure to be smooth being
that if $X$ is a compact manifold, then $\tilde{X}$ is a compact manifold
with boundary (the preimage of $Y$).  Similarly, if $X$ is a manifold with
boundary, then $\tilde{X}$ can have corners.  Recall that a smooth manifold
with corners is a (paracompact, Hausdorff) space $X$ with a covering $U_i$
by open sets, each homeomorphic to a space $\R^{p_i}\times (\R^{\ge
  0})^{q_i}$, and in such a way that the compatibility maps are $C^\infty$.
We extend this to pairs $(X,Y)$ by insisting that each $Y\cap U_i$ either
be empty or an intersection of coordinate hyperplanes.  (This is
essentially just a condition that $Y$ meets the boundary and corners of $X$
transversely.)  In particular, given a pair, we may associate a normal
bundle $N_Y(X)$, which on a patch $U_i\cap Y$ is the quotient of the cone
$\R^{p_i}\times (\R^{\ge 0})^{q_i}$ by the subspace $U_i\cap Y$.

The {\em spherical blow-up} of the pair $(X,Y)$ is then defined as follows.
If $Y$ is empty, the blow-up of $(X,Y)$ is itself; if $X$ is a cone bundle
over $Y$, then the blow-up is the pair $(\tilde{X},\tilde{Y})$, where
$\tilde{X}$ is the closure of the subset
\[
X-Y\subset X\times X/\R^+,
\]
and $\tilde{Y}\cong X/\R^+$ is the preimage of $Y$ in $\tilde{X}$.  In
general, any pair $(X,Y)$ looks locally like one of the above two examples,
and the above constructions are sufficiently compatible to give a global
construction.  In particular, note that $(\tilde{X},\tilde{Y})$ is a smooth
pair with corners, and the induced map $\tilde{X}-\tilde{Y}\to X-Y$ is a
diffeomorphism.

For instance, if $X\subset \R^n$ and $Y\subset X$ consists of a single
point, then there is a natural projection $X\setminus Y\to S^{n-1}$, and
$\tilde{X}$ is the closure of $X\setminus Y$ in $\R^n\times S^{n-1}$.  The
subset $\tilde{Y}$ can then be viewed as the set of all unit vectors
pointing towards $X$ from $Y$.

Since $\tilde{Y}\cong N_Y(X)/\R^+$, one can define other blow-ups as
push-forwards of appropriate surjections from $\tilde{Y}$.  In particular,
the real blowup corresponds to the map $N_Y(X)/\R^+\to N_Y(X)/\R^*$ (valid
so long as $N_Y(X)$ is a vector bundle, e.g., if $Y$ is disjoint from the
corners of $X$), while the complex blowup corresponds to a map
$N_Y(X)/\R^+\to N_Y(X)/\C^*$, assuming such a complex structure exists.
Both of these, unlike the spherical blow-up, preserve smoothness; however,
the spherical blow-up is a useful technical tool, both because it maps to
these cases and because it has the homotopy type of $X\setminus Y$.  (See also
\cite{GaiffiG:2004}, which considers the analogue of
${\overline Y}_{\cal G}(\R)$ replacing real blow-ups with spherical blow-ups.)

It turns out that in general, there is a long exact sequence relating the
homology of $X$, the homology of the blow-up, and the homology of the
mapping cone of the projection $\pi:\pi^{-1}(Y)\to Y$.  For a continuous
map $f:X\to Y$, recall that the mapping cone $M_f$ is defined as the
quotient of $\cone(X)\cup Y$ by identifying each point in $X$ with its
image in $Y$.  (By convention, if $X=\emptyset$, $\cone(X)$ consists of a
single point.)  Note that commutative diagrams
\[
\begin{CD}
A@>f>>B\\
@VgVV @VhVV\\
C@>k>>D
\end{CD}
\]
induce continuous maps between mapping cones; moreover, there is a natural
homeomorphism
\[
M_{(g,h):M_f\to M_k}\cong M_{(f,k):M_g\to M_h}.
\]
Given a map of pairs $f:(A,B)\to (C,D)$, define
$H_*(f):=H_*(M_{f:A\to C},M_{f:B\to D})$.  Note that if $B=D=\emptyset$,
then this is the reduced homology of $M_{f:A\to C}$; if $f$ is also an
inclusion map, then $H_*(f)\cong H_*(C,A)$.

\begin{lem}
Let $f:(A,B)\to (C,D)$ be a continuous map of pairs.  Then there is
a long exact sequence
\[
\begin{CD}
\cdots @>>> H_*(A,B) @>f_*>> H_*(C,D) @>>> H_*(f) @>>> H_{*-1}(A,B)
@>>>
\cdots
\end{CD}
\]
functorial in the sense that for any commutative square
\[
\begin{CD}
(A,B)@>f>>(C,D)\\
@VVV       @VVV\\
(E,F)@>g>>(G,H)
\end{CD},
\]
the corresponding diagram
\[
\begin{CD}
\cdots @>>> H_*(A,B) @>f_*>> H_*(C,D) @>>> H_*(f) @>>>
H_{*-1}(A,B) @>>>\cdots\\
@. @VVV @VVV @VVV @VVV @.\\
\cdots @>>> H_*(E,F) @>g_*>> H_*(G,H) @>>> H_*(g) @>>>
H_{*-1}(E,F) @>>>\cdots
\end{CD}
\]
commutes.
\end{lem}

\begin{proof}
If $B=D=\emptyset$, $A,C\ne\emptyset$, then we have the commutative diagram
\[
\begin{CD}
\cdots @. \tilde{H}_*(A) @>f_*>> \tilde{H}_*(C) @. H_*(f:A\to C)
@. \tilde{H}_{*-1}(A) @>f_*>>\cdots\\
@.   @V\sim VV @| @| @V\sim VV @.\\
\cdots @>>> H_{*+1}(M_f,C) @>>> \tilde{H}_*(C) @>>>
\tilde{H}_*(M_f)
@>>> H_{*}(M_f,C) @>>>\cdots
\end{CD}
\notag
\]
which gives rise to a reduced homology version of the desired sequence.
For the general case, let $M_{A,B}$ denote the mapping cone of the
inclusion $B\to A$, and observe that the mapping cone of
\[
f:M_{A,B}\to M_{C,D}
\]
is homeomorphic to the mapping cone of the inclusion
\[
M_{f:B\to D}\to M_{f:A\to C}.
\]
Thus we have the commutative diagram
\[
\begin{CD}
\cdots @. H_*(A,B)@>f_*>> H_*(C,D) @.  H_*(f:(A,B)\to (C,D)) @. H_{*-1}(A,B) @>f_*>> \cdots\\
@.       @V\sim VV      @V\sim VV      @V\sim VV   @V\sim VV              @.\\ 
\cdots @>>> \tilde{H}_*(M_{A,B}) @>f_*>> \tilde{H}_*(M_{C,D})
@>>> H_*(f:M_{A,B}\to M_{C,D}) @>>> \tilde{H}_{*-1}(M_{A,B}) @>>>\cdots
\end{CD}
\]
\end{proof}

\begin{rem}
More generally, given two maps
$f:(A,B)\to (C,D)$, $g:(C,D)\to (E,F)$, there is a long exact sequence
\[
\begin{CD}
\cdots @>>> H_*(f) @>(1,g)>> H_*(g\circ f) @>(f,1)>> H_*(g) @>>> H_{*-1}(f)@>>>\cdots
\end{CD}
\]
which gives the above result when $(A,B)=(\emptyset,\emptyset)$, and when
$B=D=E=F=\emptyset$ gives the reduced homology version.
\end{rem}

\begin{cor}\label{cor:l_e_s}
Suppose the map $f:(A,B)\to (C,D)$ induces an isomorphism on relative
homology.  Then there is a functorial long exact sequence
\[
\begin{CD}
\cdots @>>> H_*(A) @>f_*>> H_*(C) @>>> H_*(f:B\to D) @>>> H_{*-1}(A) @>>>
\cdots
\end{CD}
\]
\end{cor}

\begin{proof}
It suffices to show that the inclusion maps induce an isomorphism
\[
H_*(f:B\to D)\cong H_*(f:A\to C).
\]
From the long exact sequence of relative homology, this is equivalent to
the vanishing of the homology group
\[
H_*(f:(A,B)\to (C,D))
\]
which by the mapping cone long exact sequence is equivalent to the claim
that $f_*:H_*(A,B)\to H_*(C,D)$ is an isomorphism.
\end{proof}

In the case of a blow-up, we have the following.

\begin{thm}
Let $(X,Y)$ be a smooth manifold pair with corners, and let
$(\tilde{X},\tilde{Y})$ be a corresponding (spherical, real, complex)
blow-up.  Then there is a long exact sequence
\[
\cdots \to H_*(\tilde{X})\to H_*(X)\to H_*(\tilde{Y}\to Y)\to
H_{*-1}(\tilde{X})\cdots
\]
Moreover, the map $H_*(X)\to H_*(\tilde{Y}\to Y)$ factors through the
corresponding map $H_*(X)\to H_*((N_Y(X)/\R^+)\to Y)$.
\end{thm}

\begin{proof}
It suffices to show that the natural map $\pi_*:H_*(\tilde{X},\tilde{Y})\to
H_*(X,Y)$ is an isomorphism; indeed, then the long exact sequence follows
from Corollary \ref{cor:l_e_s}, and the fact that the map factors follows
from functoriality.

Now, if $U$ is an open subset of a patch, then
$\pi_*:H_*(\pi^{-1}(U),\pi^{-1}(U\cap Y))\to H_*(U,U\cap Y)$ is an
isomorphism, since $\pi^{-1}(U\cap Y)$ and $U\cap Y$ are both deformation
retracts of neighborhoods.  It then follows from Mayer-Vietoris that this
holds whenever $U$ is a finite union of such subsets, and then by taking a
direct limit, that it holds in general.
\end{proof}

The fact that the map factors through the spherical case is surprisingly
powerful, and in particular gives short exact sequences in many cases.

\begin{cor}\label{cor:blowup_connecting_trivial}
Let $(X,Y)$ and $(\tilde{X},\tilde{Y})$ be as above.  In the case of a
complex blow-up, the connecting map $H_*(X;R)\to H_*(\tilde{Y}\to Y;R)$ is
always 0.  In the case of a real blow-up, the connecting map is 0 if either
$R$ has characteristic 2 or each component of $Y$ has odd codimension, and
in general factors through multiplication by 2.
\end{cor}

\begin{proof}
Indeed, it suffices to consider the case of $\C^p$ or $\R^p$ blown up at
the origin, for which the claims are immediate.
\end{proof}

We now consider the implications of the blowup long exact sequence for real
De Concini-Procesi models.

\begin{cor}\label{cor:DCP_les}
Suppose $G\in {\cal G}$, that ${\cal G}':={\cal G}\setminus\{G\}$ is a
building set, and that $G\notin{\cal C}_{{\cal G}'}$; let
\[
d_G=\overline{Y}_{{\cal G}|_G\oplus {\cal G}/G}
\]
be the exceptional divisor, and
\[
d'_G=\overline{Y}_{{\cal G}'|_G\oplus {\cal G}/G}
\]
its image.  The homology groups of ${\overline Y}_{\cal G}(\R)$ and
${\overline Y}_{{\cal G}'}(\R)$ are related as follows.  If $G\not\subset
A$, then
\[
H_k(\overline{Y}_{\cal G})[A]
\cong
H_k(\overline{Y}_{{\cal G}'})[A],
\]
and otherwise we have the long exact sequence
\[
\begin{CD}
\cdots
@>>>
H_*(\overline{Y}_{{\cal G}|_G\oplus {\cal G}/G})[G\oplus A/G]
@>\phi_G>>
H_*(\overline{Y}_{\cal G})[A]
@>>>
H_*(\overline{Y}_{{\cal G}'})[A]
@>>>\cdots
\end{CD},
\]
where the connecting map is induced by the composition
\[
H_*(\overline{Y}_{{\cal G}'})
\to
H_*(\overline{Y}_{{\cal G}'}|d'_G)
\cong
H_*(N_{d'_G}({\cal G}')|d'_G)
\to
H_{*-1}(d_G),
\]
where $H_*(X|Y):=H_*(X,X\setminus Y)$ and the last map is induced from the
morphism $H_*(\R^k|\{0\})\cong H_{*-1}(S^{k-1})\to H_{*-1}(\P^{k-1})$.
%
%
\end{cor}

\begin{proof}
The key observation is that by Corollary \ref{cor:idems}, the map
\[
d_G=\overline{Y}_{{\cal G}|_G\oplus {\cal G}/G}
\to
d'_G=\overline{Y}_{{\cal G}'|_G\oplus {\cal G}/G}
\]
has a section, and thus the mapping cone long exact sequence breaks up into
split short exact sequences
\[
0\to H_{*+1}(d_G\to d'_G)\to H_*(d_G)\to H_*(d'_G)\to 0.
\]
We may thus identify the mapping cone homology with its image in $H_*(d_G)$:
\[
H_{k+1}(d_G\to d'_G)
\cong
\sum_{A\in {\cal G}/G}
H_k(\overline{Y}_{{\cal G}|_G\oplus {\cal G}/G})[G\oplus A].
\]
The identification of the map
\[
H_{*+1}(d_G\to d'_G)\to H_*({\overline Y}_{\cal G})
\]
follows from functoriality and the commutative diagram
\[
\begin{CD}
d_G@>>>d'_G\\
@V\phi_G VV @V\phi_G VV\\
{\overline Y}_{\cal G}@>>>{\overline Y}_{{\cal G}'}
\end{CD}.
\]

The remainder of the proof is straightforward.
\end{proof}

If $G$ is minimal in ${\cal G}$ (which is the only case that need be
considered when constructing De Concini-Procesi models via repeated
blow-ups), we find that the projective bundle is trivial, and thus we can
(noncanonically) compute the cohomology of the exceptional divisor via the
K\"unneth formula.  We find
\[
H_*(\overline{Y}_{{\cal G}|_G\oplus {\cal G}/G})[G\oplus A/G]
\cong
\Tor_\Z(\tilde{H}_{*-1}(\P_G),H_*(\overline{Y}_{{\cal G}/G})[A/G])
\oplus
\tilde{H}_*(\P_G)\otimes H_*(\overline{Y}_{{\cal G}/G})[A/G].
\]
Since the torsion subgroup of $\tilde{H}_{*-1}(\P_G)$ has exponent 2,
the $\Tor_\Z$ component consists entirely of 2-torsion; thus modulo that
2-torsion, we obtain the {\em canonical} isomorphism
\[
2 H_{*+\dim(G)-1}(\overline{Y}_{{\cal G}|_G\oplus {\cal G}/G})[G\oplus A/G]
\cong
\begin{cases}
H_{\dim(G)-1}(\P_G)\otimes 2H_*(\overline{Y}_{{\cal
    G}/G})[A/G],
& \dim(G)=0(2)\\
0
& \dim(G)=1(2).
\end{cases}
\]

\section{A cell structure for ${\overline Y}_{\cal G}$}
\label{sec:cells}

We will construct the isomorphism of the main theorem via a chain map; this
will require a careful choice of complex for the homology of ${\overline
  Y}_{\cal G}$.

To each forest $F$ in ${\cal G}$ with root $\sum_{G\in {\cal G}} G$, we
may associate a corresponding composition of operad maps:
\[
\phi_F:
\prod_i \overline{Y}_{({\cal G}|_{F_i})/\child_F(F_i)}
\to
\overline{Y}_{\cal G}
\]
Indeed, if $C$ is any sum of independent elements of $F$, we can take the
composition
\[
\phi_C\circ \phi_{F|_C}\times \phi_{F/C}.
\]
where $F|_C$ and $F/C$ are the induced forests in ${\cal G}|_C$ and ${\cal
  G}/C$.  It then follows easily from the composition law that this
composition is independent of $C$.  Since each $\phi_C$ is injective, it
follows that $\phi_F$ is injective for all $F$, and it is thus reasonable
to consider the image $d_F$ of $\phi_F$.  In particular, if $F$ consists
only of (the components of) its root, then $d_F={\overline Y}_{\cal G}$,
while if $F$ consists only of the element $G$ in addition to its root,
$d_F=d_G$.

Now, suppose ${\cal C}_{\cal G}$ has maximal element $\rt({\cal
  G})=G_1\oplus G_2\oplus \cdots\oplus G_k$ of dimension $d$.  For each
$0\le i\le d-k$ (the dimension of ${\overline Y}_{\cal G}$), let $X^i$ be
the union of $d_F$ for all forests $F$ with $\rt(F)=\rt({\cal G})$ and
$d-i$ nodes.

\begin{thm}
Suppose every chain in ${\cal C}_{\cal G}$ can be extended to a complete
flag (i.e., with all codimensions 1).  Then for all $n$,
\[
H_*(X^{n+1},X^n)
\]
is free, supported in degree $n+1$.  Moreover, the induced map
\[
H_{n+1}(X^{n+1},X^n)\to H_n(X^n)\to H_n(X^n,X^{n-1})
\]
provides these groups with a chain complex structure, with associated
homology groups canonically isomorphic to $H_*({\overline Y}_{\cal G})$.
\end{thm}

In other words, the sequence $X^0\subset X^1\subset\cdots\subset X^{d-k}$
behaves homologically as the sequence of skeletons of a CW complex.  (We
conjecture that this is in fact the sequence of skeletons of a regular cell
complex, but the homological statement will suffice for our purposes.)

\begin{proof}
It suffices to prove the claim about $H_*(X^{i+1},X^i)$, since then the
derivation of cellular homology for CW complexes carries over mutatis
mutandum.  Now, it was shown in \cite{DeConciniC/ProcesiC:1995} that the
intersection of $d_F$ and $d_{F'}$ is the submanifold $d_{F\cup F'}$ 
if $F\cup F'$ is a forest, and empty otherwise.  We thus have the canonical
isomorphism
\[
H_*(X^{n+1},X^n)
=
\bigoplus_{|F|=d-n-1} H_*(d_F,d_F\cap X^n).
\]
But then pulling this back through $\phi_F$ reduces to the case $|F|=k$;
i.e., $H_*({\overline Y}_{\cal G},\cup_{G\in {\cal G}} d_G)$.

Now, by the hypotheses, each $G$ contains a $1$-dimensional space in the
building set, or equivalently $G^\perp$ is contained in a hyperplane.  In
other words, ${\overline Y}_{\cal G}/\cup_{G\in {\cal G}}d_G$ is
homeomorphic to the quotient of $\prod_i \P_{G_i}$ by a nonempty hyperplane
arrangement.  But this in turn is homeomorphic to a wedge of spheres
(one-point compactifications of intersections of open half-spaces).
\end{proof}

\begin{rem}
Note that the hypothesis is necessary.  Indeed, suppose $F$ were a maximal
forest with $d-i$ nodes for some $i>0$ (guaranteed to exist if the
hypothesis fails).  Then $d_F\subset X^i$ is disjoint from $X^{i-1}$, and
thus its contribution to $H_*(X^i,X^{i-1})$ is simply $H_*(d_F)$.  But
$H_0(d_F)=\Z\ne 0$.
\end{rem}

With this in mind, we will denote the above chain complex as
$C_*({\overline Y}_{\cal G})$, and say that a building set ${\cal G}$
satisfying the hypotheses is ``cellular''.

For our purposes, the hypotheses are not particularly onerous; we can
simply adjoin generic hyperplanes until they hold.  This, of course, leads
to the danger that constructions based on the cellular chain complex might
depend on the choice of hyperplanes.  To control this, we may use the chain
map from the following trivial lemma.  Note that the map goes in the
reverse direction to the usual case of a morphism adjoining a subspace to a
building set.

\begin{lem}
Suppose ${\cal G}'={\cal G}\cup \{H\}$ with $\dim(H)=1$ and ${\cal G}$
cellular.  Then the identity map ${\overline Y}_{\cal G}\to {\overline
  Y}_{{\cal G}'}$ is cellular.
\end{lem}

It will also be helpful to have at our disposal various special cases of
operad maps that respect the cell structure.  The easiest is the case of a
purely operadic weak morphism.

\begin{lem}
Let ${\cal G}$ be a cellular building set, and $C\in {\cal C}_{\cal G}$.
Then ${\cal G}|_C$ and ${\cal G}/C$ are cellular, as is the operad map
$\phi_C$.  The associated map on the cellular chain complex is injective.
\end{lem}

\begin{proof}
Let $F$ be a ${\cal G}|_C\oplus {\cal G}/C$-forest with root $C\oplus
\rt({\cal G}/C)$.  Then (essentially by definition) $\phi_C\circ
\phi_F=\phi_{F'}$ for some ${\cal G}$-forest containing $C$, and this
correspondence is bijective.  In particular, $d_F$ maps isomorphically to
$d_{F'}$, and any extension of $F'$ to a $d$-node forest induces such an
extension of $F$.  The claim follows.
\end{proof}

\begin{rem}
Thus, more generally, the maps $\phi_F$ are all cellular (hardly
surprising) and act injectively on the chain complex.
\end{rem}

\begin{lem}
Let ${\cal G}'\subset {\cal G}$ be cellular building sets in $V$.  If for
all $G\in {\cal G}$, $G\cap \rt({\cal G}')\in {\cal C}_{{\cal G}'}$, then
the induced (surjective) morphism $\phi_\iota:{\overline Y}_{\cal
  G}\to{\overline Y}_{{\cal G}'}$ is cellular.  Moreover, the induced map on
the cellular chain complex is surjective.
\end{lem}

\begin{proof}
We find by a straightforward induction that for any ${\cal G}$-forest $F$
with root $\rt({\cal G})$, the image of $d_F$ under $\phi_\iota$ is of the
form $d_{F'}$, where $F'$ is the ${\cal G}'$-forest consisting of all
components of all spaces $G\cap\rt({\cal G}')$ for $G\in F$.  The induced
morphism
\[
H_{\dim(d_F)}(d_F,d_F\cap X^{\dim(d_F)-1})
\to
C_{\dim(d_F)}({\overline Y}_{{\cal G}'})
\]
is thus 0 unless $\dim(d_F)=\dim(d_{F'})$, in which case it is given
by the isomorphism
\[
H_{\dim(d_F)}(d_F,d_F\cap X^{\dim(d_F)-1})
\cong
H_{\dim(d'_F)}(d'_F,d'_F\cap X^{\dim(d'_F)-1}).
\]

It therefore remains only to show that any ${\cal G}'$-forest $F'$ with
root $\rt({\cal G}')$ can be obtained as the image of some $\phi_F$ of the
same dimension.  But this is straightforward: we may extend any forest
chain associated to $F'$ by including a complete chain from $\rt({\cal
  G}')$ to $\rt({\cal G})$, and thus obtain a forest chain with associated
forest $F$.
\end{proof}

\begin{rem}
  Suppose ${\cal G}'$ is a cellular building set, and let $G\in V$ be such
  that ${\cal G}={\cal G}'\cup \{G\}$ is a building set.  If $G\in {\cal
    C}_{{\cal G}'}$, then ${\cal G}$ is cellular, and the lemma applies.
  Otherwise, since ${\cal G}'\cup \{G\}$ is a building set, any minimal
  element of ${\cal C}_{{\cal G}'}$ containing $G$ is indecomposable, and
  thus for a generic $v\in G$, ${\cal G}'\cup \{\la v\ra\}$ and ${\cal
    G}\cup \{\la v\ra\}$ are building sets.  By induction on $\dim\rt({\cal
    C}_{{\cal G}'|_G})$, we eventually obtain a pair of building sets
  satisfying the hypotheses of the Lemma.  It follows that any map
  $\phi_\iota$ can be written as a product of morphisms each of which is
  cellular for some cell structure of the kind we are considering.
\end{rem}

Although ${\overline Y}_{\cal G}$ is not in general orientable, we can
still define a fundamental class in 
\[
H_{d-k}({\overline Y}_{\cal G},X^{d-k-1})
\]
as long as the components of the maximal element of ${\cal C}_{\cal G}$ are
even-dimensional.  Indeed, we then have a natural map
\[
\bigotimes_i \ori(G_i)
\cong
H_{d-k}(\prod_i \P_{G_i})
\to
H_{d-k}(\prod_i \P_{G_i},(\prod_i \rho_{G_i})(X^{d-k-1}))
\cong
H_{d-k}({\overline Y}_{\cal G},X^{d-k-1}).
\]
More generally, given an ordered $2$-divisible forest $F$ with
$\rt(F)=\rt({\cal G})$, we obtain a map
\[
\mu_F:
\ori(\rt({\cal G}))
\to
H_*(d_F,d_F\cap X^{d-|F|-1})
\]
by composing a map of the above form with $\phi_F$.

\begin{thm}\label{thm:calc_boundary}
Let ${\cal G}$ have root $A$, and fix a class $\omega\in \ori(A)$.
Then the map $\mu:F\mapsto \mu_F(\omega)$ gives a well-defined map
\[
\mu:C^{\dim(A)-*}_f(A)\to C_*({\overline Y}_{\cal G}),
\]
which is operadic with respect to purely operadic weak morphisms $\phi_C$
with $C$ $2$-divisible, and satisfies
\[
\partial \mu(F)
=
2\mu(dF).
\]
The corresponding map on homology is independent of the choice of
hyperplanes making ${\cal G}$ cellular.
\end{thm}

\begin{proof}
The fact that $\mu$ is operadic with respect to purely operadic weak
morphisms follows essentially immediately from its definition.
In particular, we may pull the boundary computation back through $\phi_F$,
and thus reduce to the case that $F$ consists only of root nodes.
But then the product structure of ${\overline Y}_{\cal G}$ allows us to
reduce to the case that $F$ consists of the single node $V^*$.

Since the cellular chain complex splits as a direct sum over forests,
we need only determine the induced map
\[
\begin{CD}
H_n(\P(V))
@>>>
H_n({\overline Y}_{\cal G},X^{n-1})
@>\partial >>
H_{n-1}(d_G,d_G\cap X^{n-2})
\end{CD}
\]
for $G\in {\cal G}\setminus \{V^*\}$, where $n=\dim(V)-1$.  Let $\pi$ be the
natural map ${\overline Y}_{\cal G}\to {\overline Y}_{\{G,V^*\}}$.  Then we
have a commutative diagram
\[
\begin{CD}
H_n(\P(V))
@>>>
H_n({\overline Y}_{\cal G},X^{n-1})
@>\partial>>
H_{n-1}(d_G,d_G\cap X^{n-2})
@<<< H_{n-1}(\P_G\times \P_{V^*/G})
\\
@| @V\pi VV @V\pi VV @|
\\
H_n(\P(V))
@>>>
H_n({\overline Y}_{\{G,V^*\}},\pi(X^{n-1}))
@>\partial>>
H_{n-1}(\pi(d_G),\pi(d_G\cap X^{n-2}))
@<<< H_{n-1}(\P_G\times \P_{V^*/G})
\\
@| @AAA @AAA @|
\\
H_n(\P(V))
@>>>
H_n({\overline Y}_{\{G,V^*\}},\pi(d_G))
@>\partial>>
H_{n-1}(\pi(d_G))
@= H_{n-1}(\P_G\times \P_{V^*/G})
\end{CD}
\notag
\]
in which the first row of vertical arrows consists of isomorphisms.  In
particular, it suffices to determine the action of $\partial$ in the bottom
row, which is straightforward.

That the action on homology is independent of the choice of hyperplanes
follows immediately from the fact that the maps $\mu_F$ commute with the
adjunction-of-hyperplanes maps.
\end{proof}

With this in mind, let $H^{*}_f(A;2d)$ denote the cohomology of the
forest complex with differential multiplied by 2.  Note that the only
effect of this is to extend the cohomology by some $2$-torsion:
\[
2H^*_f(A;2d)\cong H^*_f(A).
\]
Since we will be focusing on the root, it will be notationally convenient
to define $H^*_f({\cal G};2d):=H^*_f(\rt({\cal G});2d)$, and similarly for
chain complexes.

\begin{lem}
  Let $G\in {\cal G}$ be such that ${\cal G}'={\cal G}\setminus\{G\}$ is a
  cellular building set with $G\in {\cal C}_{{\cal G}'}$.  Then for any
  $2$-divisible ${\cal G}$-forest $F$, $(\phi_\iota)_*(\mu(F))$ can be
  computed as follows.  If $G\not\in F$, then
  $(\phi_\iota)_*(\mu(F))=\mu(F)$.  If $G\in F$ and there is a unique
  ${\cal G}'$-component $G'$ of $G$ which is not a component of
  $\child_F(G)$, then $(\phi_\iota)_*(\mu(F))=\mu(F')$ where $F'$ is the
  (2-divisible) forest obtained from $F$ by replacing $G$ by $G'$.  In all
  remaining cases, $(\phi_\iota)_*(\mu(F))=0$.
\end{lem}

\begin{proof}
This follows easily by a consideration of the composition
$\phi_\iota\circ\phi_F$, using the fact that
\[
\phi_\iota\circ\phi_C
=
\begin{cases}
\phi_C\circ (\phi_\iota\times 1) & G\subset C\\
\phi_C\circ (1\times \phi_\iota) & G\not\subset C;
\end{cases}
\]
moreover in the second case, the morphism $\iota$ is the identity whenever
$G+C\in {\cal G}'/C$.  The claim when $G\notin F$ is thus immediate.

For $G\in F$, we may take $C=G$ and thus reduce to considering the case
$G=\rt(F)$.  In that case, taking $C=\child_F(G)$ above, if $G =
\child_F(G)+G'$ for some $G'\in {\cal G}'$, then we may take $G'$ to be a
${\cal G}'$-component of $G$ (the unique such component not a component of
$\child_F(G)$), and since $G'/C=G/C$, we obtain $\phi_{F'}$ as desired.
Note that since $G'$ is a component of $G$ and all other components of $G$
are even-dimensional, the same applies to $G'$; in addition, since $G$ is a
sum of subspaces in $F'$, the root of $F'$ is unchanged.

Otherwise, if $G\in F$ but $G$ has multiple components not contained in
$\child_F(G)$, then we reduce to the case $F=\{G\}$, for which triviality
of the map follows by dimension considerations.
\end{proof}

\begin{lem}
Let ${\cal G}$ be a building set with root $A$.  The image in homology of
the map $\mu$ is contained in
$H_*({\overline Y}_{\cal G})[A].$
\end{lem}

\begin{proof}
  We need to show that for any $C\in {\cal C}_{\cal G}$ not equal to $A$,
  $\pi_C\circ\mu=0$.  This is a slight difficulty, given that we may not be
  able to arrange for $\pi_C$ itself to be cellular.  However, we can
  always write $\pi_C$ as a composition of maps $\phi_\iota$, each of which
  removes a single element of ${\cal G}$.  Each of those maps can be made
  cellular (and made to satisfy the hypotheses of the lemma) by adjoining
  suitably generic hyperplanes, and the action thus commutes with $\mu$.
  At some point in the chain of maps, the root of the building set
  necessarily decreases; we claim that at that point, the image of $\mu$ is
  annihilated.  But this follows from the fact that for such a map,
  $(\phi_\iota)_*(\mu(F))=0$ for all 2-divisible ${\cal G}$-forests $F$ with
  root $A$.
\end{proof}

\begin{thm}\label{thm:mu_plays_nice_with_blowup}
The map $\mu$ is compatible with the blow-up long exact sequence, as
follows.  Suppose $G\in {\cal G}$, that ${\cal G}':={\cal G}\setminus\{G\}$
is a building set, and that $G\notin{\cal C}_{{\cal G}'}$, but $G\subset
\rt({\cal G}')$.  Then we have the following commutative diagram of long
exact sequences:
\[
\begin{CD}
\cdots
@>>>
H^{\dim(A)-*}_f({\cal G}|_G\oplus {\cal G}/G;2d)
@>\phi_G>>
H^{\dim(A)-*}_f({\cal G};2d)
@>\phi_\iota>>
H^{\dim(A)-*}_f({\cal G}';2d)
@>>>
\cdots
\\
@. @V\mu VV @V\mu VV @V\mu VV @.
\\
\cdots
@>>>
H_*(\overline{Y}_{{\cal G}|_G\oplus {\cal G}/G})[G\oplus A/G]
@>\phi_G>>
H_*(\overline{Y}_{\cal G})[A]
@>\phi_\iota>>
H_*(\overline{Y}_{{\cal G}'})[A]
@>>>
\cdots
\end{CD}
\]
\notag
\end{thm}

\begin{proof}
First, by adjoining sufficient hyperplanes to ${\cal G}$, we may arrange
for the maps to be cellular; in particular, this means that now $G\in {\cal
  C}_{{\cal G}'}$, but with at least one hyperplane component.

Now, the top row is the long exact sequence associated to the short exact
sequence
\[
\begin{CD}
0
@>>>
C^{\dim(A)-*}_f({\cal G}|_G\oplus {\cal G}/G;2d)
@>\phi_G>>
C^{\dim(A)-*}_f({\cal G};2d)
@>\phi_\iota>>
C^{\dim(A)-*}_f({\cal G}';2d)
@>>>
0
\end{CD}
\notag
\]
of cochain complexes; we will thus be able to obtain the desired diagram if
we can exhibit a corresponding short exact sequence for the bottom row
making everything commute.

On the bottom row, it is tempting to consider the sequence
\[
\begin{CD}
0
@>>>
C_*(\overline{Y}_{{\cal G}|_G\oplus {\cal G}/G})
@>\phi_G>>
C_*(\overline{Y}_{\cal G})
@>\phi_\iota>>
C_*(\overline{Y}_{{\cal G}'})
@>>>
0,
\end{CD}
\notag
\]
which indeed forms a commutative diagram with the sequence on forest
cochains.  Unfortunately, this sequence is not exact in the middle; in
fact, $\phi_\iota\circ\phi_G\ne 0$.  However, we {\em do} have the
inclusion $\ker\phi_\iota\subset \im\phi_G$, and can verify that
$\phi_\iota\circ\phi_G\circ\mu=0$.  We thus obtain the commutative diagram
\[
\begin{CD}
0
@>>>
C^{\dim(A)-*}_f({\cal G}|_G\oplus {\cal G}/G;2d)
@>\phi_G>>
C^{\dim(A)-*}_f({\cal G};2d)
@>\phi_\iota>>
C^{\dim(A)-*}_f({\cal G}';2d)
@>>>
0\\
@. @V\mu VV @V\mu VV @V\mu VV @.\\
0
@>>>
\ker(\phi_\iota\circ\phi_G)
@>\phi_G>>
C_*(\overline{Y}_{\cal G})
@>\phi_\iota>>
C_*(\overline{Y}_{{\cal G}'})
@>>>
0.
\end{CD}
\notag
\]
Now, we have the composition $\phi_\iota\circ\phi_G=\phi_G\circ \pi_{A/G}$
(note that $\pi_{A/G}:{\cal G}|_G\times {\cal G}/G\to {\cal G}/G$ is
cellular) and thus $\ker(\phi_\iota\circ\phi_G)=\ker(\pi_{A/G})$.  The
claim follows.
\end{proof}

We can now prove a significant portion of the main theorem.

\begin{thm}\label{thm:main_odd}
Let ${\cal G}$ be a building set.  Then the maps
\[
\mu:
H^{\dim(C)-*}_f({\cal G}|_C;2d;\Z[1/2])\otimes \ori(C)
\to
H_*(\overline{Y}_{{\cal G}|_C};\Z[1/2])[C]
\cong
H_*(\overline{Y}_{\cal G};\Z[1/2])[C]
\]
induce an isomorphism
\[
\bigoplus_{C\in {\cal C}_{\cal G}}
H^{\dim(C)-*}_f({\cal G}|_C;2d;\Z[1/2])\otimes \ori(C)
\cong
H_*(\overline{Y}_{\cal G};\Z[1/2]).
\]
\end{thm}

\begin{proof}
It suffices to consider the case $C=\rt({\cal G})$.  If ${\cal G}$ consists
precisely of the components of $C$ (so ${\overline Y}_{\cal G}$ is a
product of projective spaces), then the isomorphism is immediate; the claim
in general follows by induction from Theorem
\ref{thm:mu_plays_nice_with_blowup} and the five-lemma.
\end{proof}

\section{$2$-adic homology}
\label{sec:Z2}

Of course, we expect (according to the main theorem) that $\mu$ should
induce an isomorphism
\[
\bigoplus_{C\in {\cal C}_{\cal G}}
2 H^{\dim(C)-*}_f({\cal G}|_C;2d;\Z)\otimes \ori(C)
\cong
2 H_*(\overline{Y}_{\cal G};\Z).
\]
Since we have shown this over $\Z[1/2]$, it remains only to consider the
case of $2$-adic coefficients.  The difficulty here is that the five-lemma
does not apply if one weakens ``isomorphism'' to ``isomorphism modulo
$2$-torsion''.  Notably, the $2$-torsion being ignored by the known maps
could easily combine to form $4$-torsion in the unknown (middle) map.  If
the middle map induces an isomorphism modulo $2$-torsion on mod 4 homology,
this cannot happen; this suggests that it should suffice for us to consider
the action of $\mu$ on mod 4 homology.

In fact, it is not necessary to use the blow-up long exact sequence to
reduce to mod 4 homology; one can simply use the following completely
general lemma.

\begin{lem}\label{lem:2adic_from_mod4}
Let $\rho:C_*\to D_*$ be a chain map of $2$-adic chain complexes.
If $\rho$ induces an isomorphism
\[
2 H_*(C\otimes \Z/4\Z)\to 2 H_*(D\otimes \Z/4\Z),
\]
then it induces an isomorphism
\[
2 H_*(C)\to 2 H_*(D).
\]
\end{lem}

\begin{proof}
Recall that there is a singly-graded spectral sequence (the Bockstein
spectral sequence) with $r$-th page
\[
B^r_n(C) = 2^{r-1} H_n(C\otimes \Z/2^r \Z)
\]
for $r\ge 1$.  The hypothesis states that $\rho$ induces an isomorphism
$B^2_*(C)\cong B^2_*(D)$; since $\rho$ is a chain map, this is an
isomorphism of spectral sequences, and thus $\rho$ induces an isomorphism
$B^r_*(C)\cong B^r_*(D)$ for all $r\ge 2$.

Now, since the short exact sequence from the universal coefficient theorem
is split, it follows that it remains exact when multiplied by a power of 2.
In particular, we have a (split) short exact sequence
\[
0
\to
2^{r-1}(H_n(C)\otimes \Z/2^r \Z)
\to
2^{r-1} H_n(C\otimes \Z/2^r\Z)
\to
2^{r-1} \Tor(H_{n-1}(C),\Z/2^r\Z)
\to
0.
\]

Now, suppose that we know that $\rho$ induces an isomorphism $2H_{n-1}(C)\cong
2H_{n-1}(D)$ (certainly the case for $n\le 0$).  Then it follows from the
universal coefficient exact sequence that it induces isomorphisms
\[
2^{r-2}((2 H_n(C))\otimes \Z/2^{r-1} \Z)
=
2^{r-1}(H_n(C)\otimes \Z/2^r \Z)
\cong
2^{r-1}(H_n(D)\otimes \Z/2^r \Z)
=
2^{r-2}((2 H_n(D))\otimes \Z/2^{r-1} \Z)
\]
for all $r\ge 2$.  But then, by the five lemma, $\rho$ induces an
isomorphism $2H_n(C)\cong 2H_n(D)$ as required, and the result follows by
induction.
\end{proof}

Now, $2(H_*(C\otimes \Z/4\Z))$ can be computed as follows: Let
\[
\beta:H_*(C\otimes \F_2)\to H_{*-1}(C\otimes \F_2)
\]
be the Bockstein morphism, i.e., the connecting map in the long exact
sequence
\[
\begin{CD}
\cdots
@>>>
H_*(C\otimes \Z/4\Z)
@>>>
H_*(C\otimes \F_2)@>\beta>> H_{*-1}(C\otimes \F_2)@>>> H_{*-1}(C\otimes
\Z/4\Z)
@>>>
\cdots
\end{CD}
\]
Then $\beta$ makes $H_*(C\otimes \F_2)$ a chain complex, and
\[
2 H_*(C\otimes \Z/4\Z) \cong H_*(H_*(C\otimes \F_2);\beta).
\]
In particular, we may restate the lemma by saying that $\rho$ induces an
isomorphism modulo 2-torsion on homology iff it induces an isomorphism on
$\beta$-homology.

We thus need to understand the mod 2 homology of ${\overline Y}_{\cal G}$,
and the associated Bockstein morphism.  Now, since ${\overline Y}_{\cal G}$
is a smooth real algebraic variety, every smooth subvariety induces a class
in mod 2 homology (since issues of orientability can be ignored).  In
general, such classes do not span homology, but in our case, they do
suffice; i.e., ${\overline Y}_{\cal G}$ is ``algebraically maximal'' in the
sense of \cite{KrasnovVA:2003}.  This follows, for instance, from the fact
\cite{DeConciniC/ProcesiC:1995} that the cohomology of ${\overline Y}_{\cal
  G}(\C)$ is generated by $\R$-rational cycles.  (This can also be seen
directly from Proposition \ref{prop:mod2basis} below, which gives an
explicit basis of classes associated to subvarieties.)

This allows us to extend the construction of the fundamental class to
classes of other degrees.  If $V^*\in {\cal G}$ and $\dim(V)=d$, we
construct a canonical class in $H_i({\overline Y}_{\cal G};\F_2)$ as
follows.  Choose a generic $i+1$-dimensional subspace of $V^*$, and
consider the closure of its image in ${\overline Y}_{\cal G}$.  The result
is an $i$-dimensional subvariety, and any two such subvarieties are clearly
homotopic; we thus obtain the desired canonical homology class.  Moreover,
we readily verify that if $i>0$, this class is actually in $H_i({\overline
  Y}_{\cal G};\F_2)[V^*]$.

More generally, let $F$ be a forest (not necessarily $2$-divisible) with
root $\rt({\cal G})$, and suppose we are also given a map $d:F\to \Z^+$
such that $d(G)<\dim(G)-\dim(\child_F(G))$ for all $G\in F$.  Then we have
a canonical class in
\[
H_{d(G)}(({\cal G}|_G)/\child_F(G);\F_2)
\]
for all $G\in F$, and may define $\mu(F,d)$ to be the image under $\phi_F$
of the product of these classes.  Note that this is compatible with our
previous notation, in the sense that if $F$ is $2$-divisible and $d_F(G) =
\dim(G)-\dim(\child_F(G))-1$, then $\mu(F,d_F)=\mu(F)\otimes \F_2$.

\begin{prop}\label{prop:mod2basis}
The classes $\mu(F,d)$ for $\kappa(F,d):=\sum_{G\in F} d(G)=k$ form a basis
of $H_k({\overline Y}_{\cal G};\F_2)[\rt({\cal G})]$.
\end{prop}

\begin{proof}
Let $A:=\rt({\cal G})$.  Modulo 2, the blow-up exact sequence decomposes
into short exact sequences
\[
\begin{CD}
0
@>>>
H_*(\overline{Y}_{{\cal G}|_G\oplus {\cal G}/G};\F_2)[G\oplus A/G]
@>\phi_G>>
H_*(\overline{Y}_{\cal G};\F_2)[A]
@>>>
H_*(\overline{Y}_{{\cal G}'};\F_2)[A]
@>>>
0
\end{CD}
\]
where we recall that $G\in {\cal G}$ such that ${\cal G}':={\cal
  G}\setminus\{G\}$ is a building set and $G\notin {\cal C}_{{\cal G}'}$;
the connecting map is trivial by Corollary
\ref{cor:blowup_connecting_trivial}.  We moreover readily verify that if
$G\in F$, so there is a class
\[
\mu(F,d)
\in
H_*(\overline{Y}_{{\cal G}|_G\oplus {\cal G}/G})[G\oplus A/G],
\]
then $\phi_G(\mu(F,d))=\mu(F,d)$.  Similarly, if $G\notin F$, then
$\mu(F,d)$ maps to its counterpart in $H_*(\overline{Y}_{{\cal G}'})[A]$.
But then the claim follows by induction.
\end{proof}

\begin{rem}
  The above basis (summed over the graded pieces) is trivially bijective
  with the basis of $H^*({\overline Y}_{\cal G}(\C))$ given in
  \cite{YuzvinskyS:1997}, which by the results of \cite{KrasnovVA:2003}
  alluded to above induces a basis of $H_*({\overline Y}_{\cal G};\F_2)$.
  However, the precise relation between these bases is somewhat unclear;
  for instance, Yuzvinsky's basis is not compatible with the natural
  grading by ${\cal C}_{\cal G}$.  Note also that precisely the same
  argument (the long exact sequence for complex blow-ups also splits into
  short exact sequences) shows that the above algebraic cycles form a basis
  of the homology of the complex locus.
\end{rem}

Next, we need to determine the action of the Bockstein morphism.
For $G\in {\cal G}$, we define a map $\beta_G$ by
\[
\beta_G(\mu(F,d))=0
\]
if $G\in F$, $F\cup\{G\}$ is not a forest, or
\[
\dim(G/\child_{F\cup\{G\}}(G))\bmod 2=1.
\]
Otherwise,
\[
\beta_G(\mu(F,d))=\mu(F\cup\{G\},d')
\]
where
\[
d'(H)=\begin{cases}
\dim(G/\child_{F\cup\{G\}}(G))-1 & H=G\\
d(H)-\dim(G/\child_{F\cup\{G\}}(G)) & G\in \child{H}\\
d(H) & \text{otherwise}.
\end{cases}
\]

\begin{lem}
The Bockstein morphism acts on forest classes $\mu(F,d)$ as follows:
\[
\beta(\mu(F,d))
=
\sum_{G\in {\cal G}} \beta_G(\mu(F,d))
+
\sum_{\substack{G\in F\\d(G)\bmod 2=0}} \mu(F,d-\delta_G).
\]
\end{lem}

\begin{proof}
First suppose that $F=\{V^*\}$ and $d(V^*)=\dim(V)-1$.
If $\dim(V)$ is even, then this lifts to the integral chain $\mu(\{V^*\})$,
and the claim follows from Theorem \ref{thm:calc_boundary};
a similar calculation establishes the claim when $\dim(V)$ is odd.

Now, if $F=\{V^*\}$, $0<d(V^*)<\dim(V)-1$, then $\mu(F,d)$ can be
calculated using the fact that $\beta$ is functorial, so can be transported
through operad maps; in particular, the operad map associated to inclusion
of a generic $d(V^*)$-plane.  The claim for general $F$ then follows using
the map $\phi_F$.
\end{proof}

Given a pair $(F,d)$, define the {\em defect}
\[
\delta(F,d)=|\{G:G\in F|d(G)\bmod 2=0\}|.
\]
Of the two components of the above expression for $\beta$, only the second
component changes the defect, decreasing it by 1.  We may thus interpret
$\beta$ as the differential of the total complex of a suitable double
complex.  To be precise, for $0\le p\le q$, define
\begin{align}
E^0_{p,q} &= \bigoplus_{\substack{\delta(F,d)=p\\ \kappa(F,d)=p+q}}
\F_2\mu(F,d)\\
\partial_1(\mu(F,d)) &= \sum_{\substack{G\in F\\d(G)\bmod 2=0}}
\mu(F,d-\delta_G)\\
\partial_2(\mu(F,d)) &= \sum_{G\in {\cal G}} \beta_G(\mu(F,d))
\end{align}
Then $\partial_1(E^0_{p,q})\subset E^0_{p-1,q}$,
$\partial_2(E^0_{p,q})\subset E^0_{p,q-1}$, and we readily verify that
this defines a double complex with total complex $(H_k,\beta)$.

\begin{lem}
The map
\[
\mu:
H^{\dim(C)-*}_f({\cal G};2d;\F_2)
\to
H_*(\overline{Y}_{\cal G};\F_2)[\rt({\cal G})]
\]
induces an isomorphism on Bockstein homology.
\end{lem}

\begin{proof}
  Since $\mu$ is injective and the image of the chain map $\mu$ is
  annihilated by $\partial_1$ and closed under $\partial_2$, it suffices to
  show that the classes $\mu(F)\in E^0_{0,q}$ for $2$-divisible forests $F$
  are representatives for the $\partial_1$-homology on $E^0_{p,q}$.

Now, the action of $\partial_1$ leaves the forest unchanged, so we may
restrict our attention to a single forest $F$; in other words, it suffices
to compute the homology of $\partial_1$ on the space $\bigoplus_d
\mu(F,d)$.  But this new complex is a product complex, with one factor for
each element of $G$.  For each $G$, if $\dim(G)-\dim(\child(G))$ is odd,
the homology of the corresponding factor is trivial, and thus $F$ must be
$2$-divisible.  Similarly, if $\dim(G)-\dim(\child(G))$ is even, the
homology is supported in the top degree; the lemma follows.
\end{proof}

By Lemma \ref{lem:2adic_from_mod4}, this implies that $\mu$
induces an isomorphism modulo $2$-torsion on $2$-adic homology, which
together with Theorem \ref{thm:main_odd} implies that $\mu$ induces an
isomorphism modulo $2$-torsion on integral homology.

\begin{thm}\label{thm:main_isom}
Let ${\cal G}$ be a building set.  Then the maps
\[
\mu:
H^{\dim(C)-*}_f({\cal G}|_C;2d;\Z)\otimes \ori(C)
\to
H_*(\overline{Y}_{{\cal G}|_C};\Z)[C]
\cong
H_*(\overline{Y}_{\cal G};\Z)[C]
\]
induce an isomorphism
\[
\bigoplus_{C\in {\cal C}_{\cal G}}
2H^{\dim(C)-*}_f({\cal G}|_C;2d;\Z)\otimes \ori(C)
\cong
2H_*(\overline{Y}_{\cal G};\Z).
\]
\end{thm}

\section{Operadicity}\label{sec:operadic}

To finish proving the main theorem, it remains only to prove that the
isomorphism given by $\mu$ is operadic.  There is a slight subtlety here,
in that $\mu$ is only operadic once we have quotiented out the $2$-torsion.

There are, however, three important special cases of operadicity that hold
directly on integral homology.  First, if $f:V\to W$ is a morphism such
that $f^*$ induces a bijection between ${\cal G}$ and ${\cal G}'$, then
$\phi_f$ is an isomorphism, and commutativity is trivial.  Next, if ${\cal
  G}={\cal G}'\cup \{G\}$, then we have a
commutative diagram
\[
\begin{CD}
H^{\dim(A)-*}_f({\cal G};2d)
@>\phi_\iota>>
H^{\dim(A)-*}_f({\cal G}';2d)
\\
@V\mu VV @V\mu VV
\\
H_*(\overline{Y}_{\cal G})[A]
@>\phi_\iota>>
H_*(\overline{Y}_{{\cal G}'})[A]
\end{CD}
,
\]
when $A=\rt({\cal G})=\rt({\cal G}')$, which immediately implies
commutativity of the diagram
\[
\begin{CD}
H^{\dim(A)-*}_f({\cal G}|_A;2d)
@>\phi_\iota>>
H^{\dim(A)-*}_f({\cal G}'|_A;2d)
\\
@V\mu VV @V\mu VV
\\
H_*(\overline{Y}_{\cal G})[A]
@>\phi_\iota>>
H_*(\overline{Y}_{{\cal G}'})[A]
\end{CD}
\]
for general $A\in {\cal C}_{{\cal G}'}$; together with the previous case
(and the fact that the maps are trivial if $A\notin {\cal C}_{{\cal G}'}$),
this shows operadicity whenever $f^*:{\cal G}'\to {\cal G}$ is injective.
Finally, if $C\in {\cal C}_{\cal G}$, then we have commutativity of
\[
\begin{CD}
H^{\dim(A)-*}_f({\cal G}|_C\oplus ({\cal G}/C)|_{A/C};2d)
@>\phi_C>>
H^{\dim(A)-*}_f({\cal G}|_A;2d)
\\
@V\mu VV @V\mu VV
\\
H_*(\overline{Y}_{{\cal G}|_C\oplus {\cal G}/C})[C\oplus A/C]
@>\phi_C>>
H_*(\overline{Y}_{\cal G})[A]
\end{CD}
\]
for all $A\supset C$.  Composing these three cases, we find that the diagram
\[
\begin{CD}
H^{\dim(\ker(f^*))+\dim(A)-*}_f({\cal G}'|_{\ker(f^*)}\oplus {\cal G}|_A;2d)
@>\phi_f>>
H^{\dim(C)-*}_f({\cal G}'|_C;2d)
\\
@V\mu VV @V\mu VV
\\
H_*(\overline{Y}_{{\cal G}'|_{\ker(f^*)}\oplus {\cal G}/C})[\ker(f^*)\oplus A]
@>\phi_f>>
H_*(\overline{Y}_{\cal G})[C]
\end{CD}
\]
is commutative whenever $\ker(f^*)\in {\cal C}_{{\cal G}'|_C}$; in
particular, if $A\ne f^*(C)$, the horizontal maps are 0.

Now, in general, given an operad map $f:{\cal G}\to {\cal G}'$
and spaces $A\in \Pi^{(2)}_{\cal G}$, $B,C\in \Pi^{(2)}_{{\cal G}'}$
with $B\subset \ker(f^*)$, we need to show commutativity of the diagram
\[
\begin{CD}
2H^{\dim(A)-*}_f({\cal G}'|_B\oplus {\cal G}|_A;2d)
@>\phi_f>>
2H^{\dim(A)-*}_f({\cal G}'|_C;2d)
\\
@V\mu VV @V\mu VV\\
2 H_*(\overline{Y}_{{\cal G}'|_{\ker(f^*)}\oplus {\cal G}})[B\oplus A]
@>\phi_f>>
2 H_*(\overline{Y}_{{\cal G}'})[C]
\end{CD}
\]
Since $B\subset \ker(f^*)$, we may write $f$ as a composition
\[
\begin{CD}
f:{\cal G}@>g>> {\cal G}'/B@>i_B>> {\cal G}'
\end{CD}
\]
and thus obtain a commutative diagram
\[
\begin{CD}
2H_*({\overline Y}_{{\cal G}'|_B\oplus({\cal G}'/B)|_{\ker(g^*)}\oplus {\cal G}})[B\oplus 0\oplus A]
@>{1\times \phi_g}>>
2H_*({\overline Y}_{{\cal G}'|_B\oplus ({\cal G}'/B)})[B\oplus C/B]
\\
@V{\phi_B\times 1}VV  @VV{\phi_B}V
\\
2 H_*({\overline Y}_{{\cal G}'|_{\ker{f^*}}\oplus {\cal G}})[B\oplus A]
@>{\phi_f}>>
2 H_*({\overline Y}_{{\cal G}'})[C]
\end{CD}
\]
of operad maps (the grading on the upper-right-hand corner being forced by
the requirement that $\phi_B$ act nontrivially).  Now, the vertical arrows
are respected by $\mu$ (and the left vertical arrow is an isomorphism), so
it suffices to show that $\mu$ respects the top row.  Since $\mu$ respects
products, we find that we may restrict our attention to the case $B=0$.
Similarly, we may compose with the map $\phi_C$ to restrict to the case
$C=\rt({\cal G}')$.  If $\ker(f^*)=0$, then $f^*:{\cal G}'\to {\cal G}$ is
injective, and we have already verified operadicity.  Otherwise, let
$H\subset \ker(f^*)$ be a hyperplane, which we may remove from ${\cal G}'$
if necessary to insure that the corresponding inclusion map is a morphism.
Thus if we factor
\[
\begin{CD}
f:{\cal G}@>g>> {\cal G}'/H@>i_H>> {\cal G}',
\end{CD}
\]
the composition law reads
\[
\phi_{i_H}\circ \phi_g = \phi_f\circ (\phi_{i_H|_{\im(g)}}\times 1).
\]
Since any morphism induces an isomorphism on $0$-graded homology, $\phi_f$
respects $\mu$ if $\phi_{i_H}$ and $\phi_g$ both do so.

By induction, it thus remains only to consider operadicity of the morphism
$\phi_{i_H}$.  At the poset level, this map is trivial, so we need simply
show that the map
\[
\begin{CD}
2 H_*({\overline Y}_{{\cal G}'/H})[A]
@>{\phi_{i_H}}>>
2 H_*({\overline Y}_{{\cal G}'})[C]
\end{CD}
\]
is trivial whenever $H\subset C$.  If $A\ne C/H$, this map is trivial even
before quotienting by the 2-torsion; if $A=C/H$, then one of the two spaces
has odd dimension, and thus by Theorem \ref{thm:main_isom} that graded
piece of the integral homology consists entirely of $2$-torsion, so the
map is again necessarily trivial.

It follows that $\mu$ is operadic modulo 2-torsion in general, concluding
the proof of Theorem \ref{thm:main}.

\section{Further directions}\label{sec:further}

In \cite{GaiffiG:2004}, Gaiffi considers the analogue of De Concini-Procesi
models in which the projective spaces of the definition are replaced by
spheres, and shows that the result is a manifold with corners which is
homotopy equivalent to the complement of the subspace arrangement.  Much of
the above structure carries over to this case, with two important
exceptions.  First, although the resulting structures indeed form a
(universal) topological operad, Corollary \ref{cor:idems} does not hold.
To be precise, if $C$ has codimension 1 in $\rt({\cal G})$, then the operad
map gives rise to two splittings of $\pi_C$ which are not homotopic to each
other.  Thus, in order to obtain a natural grading on homology, it is
necessary to assume that this codimension 1 situation never arises.  Of
course, this condition already appeared in
\cite{GoreskyM/MacPhersonR:1988,YuzvinskyS:2002,DeligneP/GoreskyM/MacPhersonR:2000,deLonguevilleM/SchultzCA:2001},
for similar reasons.

The second, more serious, difficulty is that since we can no longer add
hyperplanes without changing the geometry, we do not have a good analogue
of the cell structure of Section \ref{sec:cells}.  One way to avoid this
would be to generalize the problem by allowing a combination of both
projective spaces and spheres in constructing the compactification; at that
point, one could again add hyperplanes, so long as they were associated to
projective spaces.  Presumably, some constraints would be necessary on this
association in order to define the operad maps; one likely candidate is
that if $\rho_G$ maps to a projective space, then so does $\rho_H$ for all
$H\subset G$.  Also, the argument of Section \ref{sec:operadic} would no
longer apply, so it might be more difficult to prove that the associated
maps were operadic.  The relevant poset in this case would presumably be
the subposet of ${\cal C}_{\cal G}$ consisting of spaces each component of
which is associated to an orientable target ($S(G)$ or $\P(G)$ with
$\dim(G)$ even).

One important special case of this would be when $V^*\in G$, and only
$\rho_{V^*}$ maps to a sphere; this is a (nonorientable) double cover of
${\overline Y}_{\cal G}$ studied for $A_{n-1}$ in \cite{KapranovMM:1993} and
for general Coxeter arrangements in
\cite{ArmstrongSM/CarrM/DevadossSL/EnglerE/LeiningerA/ManapatM:2009}.  This
is essentially equivalent to studying the homology of the pullback of the
orientation sheaf on $\P(V)$, which is straightforward enough using the
present methods.  We obtain the following.

\begin{thm}\label{thm:twist}
Let ${\cal G}$ be a building set containing $V^*$, and let $\omega$ be the
pullback through $\rho_V$ of the orientation sheaf on $\P(V)$.  Then there
is an isomorphism
\[
2 H_*({\overline Y}_{\cal G};\omega)
\cong
\begin{cases}
0, & \dim(V)\bmod{2}=0\\
H^{\dim(V)-*}_{(1)}(\Pi^{(2)}_{\cal G}\cup \{V^*\})\otimes \ori(V),
&\dim(V)\bmod{2}=1
\end{cases}
\]
In particular, if ${\tilde Y}_{\cal G}$ is the double cover over
${\overline Y}_{\cal G}$ associated to $\omega$, then
\[
H_*({\tilde Y}_{\cal G};\Z[1/2])
\cong
H_*({\overline Y}_{\cal G};\omega\otimes \Z[1/2])
\oplus
H_*({\overline Y}_{\cal G};\Z[1/2])
\]
\end{thm}

\begin{proof} (Sketch)
The proof of the first claim simply follows the proof of Theorem
\ref{thm:main}, twisting by $\omega$ where appropriate.  For the second
claim, we observe that the short exact sequence
\[
0
\to
H_*({\overline Y}_{\cal G};\Z[1/2])
\to
H_*({\tilde Y}_{\cal G};\Z[1/2])
\to
H_*({\overline Y}_{\cal G};\omega\otimes \Z[1/2])
\to
0
\]
respects the grading by ${\cal C}_{\cal G}$, and in each graded component,
either $H_*({\overline Y}_{\cal G};\Z[1/2])[A]$ or $H_*({\overline Y}_{\cal
  G};\omega\otimes \Z[1/2])[A]$ or both must vanish.
\end{proof}

\medskip

Another direction of potential generalization follows from the observation
that the homology (modulo 2-torsion) depends only on the combinatorial
data, namely the poset $\Pi^{(2)}_{\cal G}$.  This suggests that there
should be a more combinatorial construction of ${\overline Y}_{\cal G}$
itself.  More precisely, it should be possible to generalize the main
theorem (or even the mixed generalization considered above) to the case of
a building set (of flats) in an oriented matroid (a combinatorial
generalization of a real hyperplane arrangement, see
\cite{BjornerA/LasVergnasM/SturmfelsB/WhiteN/ZieglerGM:1999}).  In
particular, the cell structure considered above seems particularly amenable
to a description in oriented matroid terms, which would ideally give a
proof that it actually corresponds to a CW complex.  (Note in particular
the fact that the proof of the Topological Representation Theorem (see,
e.g.,
\cite[Thm. 5.2.1]{BjornerA/LasVergnasM/SturmfelsB/WhiteN/ZieglerGM:1999})
attaches to any oriented matroid a natural regular cell decomposition of
the sphere such that the antipode map is cellular.  It should thus be
straightforward to use this to construct analogues of De Concini-Procesi
models in such a way as to agree with the geometric construction when the
oriented matroid comes from an actual hyperplane arrangement.)

\medskip

Finally, from a more algebraic perspective, it is worth noting that the
condition that ${\cal G}$ be real is certainly not necessary for
${\overline Y}_{\cal G}$ to be defined over $\R$ and thus give rise to a
smooth $\R$-manifold.  Indeed, all that is truly necessary is that ${\cal
  G}$ be closed under complex conjugation.  This has the effect of
replacing certain pairs of real blow-ups by complex blow-ups, but the
corresponding long exact sequences are still quite reasonable.  The main
difficulty is, once again, that the cell structure does not carry over,
making it difficult to define the appropriate homology classes via a chain
map.  In addition, the action of the operad maps on homology is much more
subtle, and in particular, need not respect the grading (even in the
special case of ${\overline Y}_{\cal G}(\C)$ viewed as a real manifold).

We make the following conjecture in this case.  Given an $\R$-rational
building set ${\cal G}$, we say that a real subspace $C\in {\cal C}_{\cal
  G}$ is {\em purely complex} if its decomposition is of the form
\[
C = \bigoplus_i G_i\oplus \overline{G_i}.
\]
Similarly, we denote by $\R({\cal G})$ the subset of ${\cal G}$ consisting
of real subspaces, and note that this is again a building set. Finally, for
real $A\in {\cal C}_{\cal G}$, let $\C(A)$ denote the sum of the complex
components of $A$.

\begin{conj}\label{conj:complex}
  Let ${\cal G}$ be an $\R$-rational building set.  Then $H_*({\overline
    Y}_{\cal G}(\R),\Z)$ is naturally graded by $\R({\cal
    C}_{\cal G})$, and there is an isomorphism
\[
2 H_*({\overline Y}_{\cal G}(\R),\Z)[A]
\cong
\!\!\!\bigoplus_{\substack{\C(A)\subset C\subset A\\C\text{ purely complex}}}\!\!\!
H_*({\overline Y}_{{\cal G}|_C}(\R),\Z)[C]
\otimes
2 H_*({\overline Y}_{\R({\cal G}/C)}(\R),\Z)[A/C].
\]
\end{conj}

\begin{rem}
  Note that ${\overline Y}_{{\cal G}|_C}$ is (as a real algebraic variety)
  the restriction of scalars of a complex De Concini-Procesi model, so
  topologically is homeomorphic to that complex variety.  In particular,
  the homology of its real locus is free; see \cite{YuzvinskyS:1997} or the
  remark following Proposition \ref{prop:mod2basis} above for an explicit
  description of a basis.
\end{rem}

\end{document}